\RequirePackage{fix-cm} 
\documentclass[a4paper, twoside, reqno, 12pt]{amsart}
\usepackage{fixltx2e}   

\usepackage[latin1]{inputenc}
\usepackage[T1]{fontenc}

\usepackage{eucal}
\usepackage{esint}
\usepackage{dsfont}
\usepackage{xspace}
\usepackage{amsgen}
\usepackage{amsthm}
\usepackage{amssymb}
\usepackage{amsmath}
\usepackage{upgreek}
\usepackage{amsfonts}
\usepackage{MnSymbol}
\usepackage{mathrsfs}
\usepackage{mathtools}
\usepackage[nice]{nicefrac}

\usepackage{a4wide}
\usepackage{changepage}

\usepackage[dvipsnames,table]{xcolor}

\definecolor{oneblue}{rgb}{0.0, 0.0, 0.85}
\definecolor{bluepigment}{rgb}{0.2, 0.2, 0.6}
\definecolor{darkgrey}{rgb}{0.273, 0.281, 0.30}
\definecolor{Lightgray}{rgb}{0.89, 0.89, 0.89}
\definecolor{Lightblue}{RGB}{214, 214, 214}
\definecolor{bckg}{RGB}{20.8, 20.8, 20.8} 
\definecolor{charcoal}{rgb}{0.21, 0.27, 0.31}
\definecolor{darkelectricblue}{rgb}{0.33, 0.41, 0.47}

\usepackage[sort&compress, comma, square, numbers]{natbib}

\usepackage{psfrag}
\usepackage{graphicx}
\usepackage{subfigure}
\usepackage{morefloats}
\usepackage{indentfirst}

\usepackage{acronym}
\usepackage{microtype}
\usepackage[labelsep=period,%
            labelfont={bf,sf,color=bluepigment},%
            justification=raggedright]{caption}

\usepackage[perpage, symbol]{footmisc}

\usepackage{pstricks}
\usepackage{epsfig}
\usepackage{pst-grad} 
\usepackage{pst-plot} 

\usepackage[colorlinks,
           urlcolor=oneblue,
           linkcolor=oneblue,
           citecolor=NavyBlue,
           bookmarksopen=false,
           pdfpagemode=UseNone,
           pagebackref]{hyperref}

\usepackage[explicit]{titlesec}

\titleformat{\section}
  {\color{NavyBlue}\Large\sffamily\bfseries}
  {}
  {0em}
  {\colorbox{bckg!5}{\parbox{\dimexpr\linewidth-2\fboxsep\relax}{\centering\thesection. #1}}}
  [\vspace*{0.33em}]

\titleformat{name=\section,numberless}
  {\color{NavyBlue}\Large\sffamily\bfseries}
  {}
  {0.0em}
  {\colorbox{bckg!10}{\parbox{\dimexpr\linewidth-2\fboxsep\relax}{\centering#1}}}
  [\vspace*{0.33em}]

\titleformat{\subsection}
  {\color{NavyBlue}\large\sffamily\bfseries}
  {}
  {0.0em}
  {\colorbox{bckg!5}{\parbox{\dimexpr\linewidth-2\fboxsep\relax}{\centering\thesubsection. #1}}}
  [\vspace*{0.33em}]

\titleformat{name=\subsection,numberless}
  {\color{NavyBlue}\Large\sffamily\bfseries}
  {}
  {0em}
  {\colorbox{bckg!10}{\parbox{\dimexpr\linewidth-2\fboxsep\relax}{\centering#1}}}
  [\vspace*{0.33em}]

\titleformat{\subsubsection}
  {\color{bluepigment}\sffamily\normalsize\bfseries}
  {\thesubsubsection}
  {0.5em}
  {#1}
  [\vspace*{0.33em}]

\titleformat{\paragraph}[runin]
  {\color{bluepigment}\sffamily\small\bfseries}
  {}
  {0em}
  {#1}

\titlespacing{\section}{1.0em}{1.5em plus 2pt minus 2pt}{1.0em plus 2pt minus 2pt}[0em]
\titlespacing{\subsection}{1.0em}{1.5em plus 2pt minus 2pt}{1.0em}[0em]
\titlespacing{\subsubsection}{1.0em}{1.5em plus 2pt minus 2pt}{1.0em plus 2pt minus 2pt}[0em]

\usepackage{titletoc}

\contentsmargin{0.5em}
\newlength{\tocsep} 
\titlecontents{section}[\tocsep]
  {\addvspace{10pt}\bfseries\sffamily}
  {\contentslabel[\thecontentslabel]{\tocsep}}
  {}
  {\ \titlerule*[0.75pc]{.}\ \thecontentspage}
  []
\titlecontents{subsection}[\tocsep]
  {\addvspace{8pt}\sffamily}
  {\contentslabel[\thecontentslabel]{\tocsep}}
  {}
  {\ \titlerule*[0.5pc]{.}\ \thecontentspage}
  []
\titlecontents*{subsubsection}[\tocsep]
  {\addvspace{2pt}\footnotesize\sffamily}
  {}
  {}
  {\ \titlerule*[0.35pc]{.}\ \thecontentspage}
  [\\*]

\makeatletter
\def\@setauthors{%
  \begingroup
  \def\thanks{\protect\thanks@warning}%
  \trivlist
  \centering\footnotesize \@topsep30\p@\relax
  \advance\@topsep by -\baselineskip
  \item\relax
  \author@andify\authors
  \def\\{\protect\linebreak}%
  \textsc{\normalsize\textcolor{darkelectricblue}{\authors}}%
  \ifx\@empty\contribs
  \else
    ,\penalty-3 \space \@setcontribs
    \@closetoccontribs
  \fi
  \endtrivlist
  \endgroup
}
\def\@settitle{\begin{center}%
  \baselineskip14\p@\relax
    \bfseries
    \textsc{\Large\textcolor{charcoal}{\@title}}
  \end{center}%
}
\makeatother

\usepackage{enumitem}
\setlist[description]{%
  topsep=30pt,               
  itemsep=5pt,               
  font={\bfseries\sffamily\color{bluepigment}}, 
}

\usepackage{fancyhdr}
\usepackage{lastpage}

\newcommand*\Title{\textcolor{bluepigment}{Conservation laws with moving nodes}}
\newcommand*\Authors{\textcolor{bluepigment}{G.~Khakimzyanov, D.~Dutykh \etal}}
\newcommand*{\plogo}{\textcolor{gray}{{\texttt{arXiv.org} / \textsc{hal}}}} 

\newcommand\invisiblesection[1]{%
  \addcontentsline{toc}{section}{#1}%
  \sectionmark{#1}}

\numberwithin{equation}{section}
\newtheorem{remark}{Remark}
\newtheorem{theorem}{Theorem}

\newcommand{\up}[1]{$^{\mathrm{\small\textsf{#1}}}$} 

\newcommand{\ab}{\bar{a}}
\newcommand{\wb}{\bar{w}}
\newcommand{\fh}{\hat{f}}

\newcommand{\gt}{\tilde{g}}
\newcommand{\G}{\mathbf{G}}

\newcommand{\R}{\mathbb{R}}
\newcommand{\Z}{\mathbb{Z}}

\newcommand{\Id}{\mathbb{I}}
\newcommand{\A}{\mathscr{A}}
\newcommand{\D}{\mathcal{D}}
\newcommand{\I}{\mathcal{I}}

\newcommand{\ue}{\mathrm{e}}

\newcommand{\Ll}{\mathscr{L}}

\renewcommand{\u}{\mathbf{u}}

\newcommand{\Dd}{\mathscr{D}}
\newcommand{\Pp}{\mathscr{P}}
\newcommand{\Rr}{\mathscr{R}}
\renewcommand{\leq}{\leqslant}
\renewcommand{\geq}{\geqslant}
\newcommand{\eps}{\varepsilon}
\renewcommand{\O}{\mathcal{O}}

\renewcommand{\S}{\mathcal{S}}
\newcommand{\f}{\boldsymbol{f}}

\renewcommand{\v}{\boldsymbol{v}}
\newcommand{\const}{\mathrm{const}}
\newcommand{\No}{$\mathrm{N}^\circ$}
\newcommand{\C}{\boldsymbol{\mathrm{C}}}
\renewcommand{\gamma}{\boldsymbol{\upgamma}}

\newcommand{\cf}{\emph{cf.}\xspace}
\newcommand{\vs}{\emph{vs.}\xspace}
\newcommand{\ie}{\emph{i.e.}\xspace}
\newcommand{\eg}{\emph{e.g.}\xspace}
\newcommand{\etc}{\emph{etc.}\xspace}
\newcommand{\etal}{\emph{et al.}\xspace}

\newcommand{\sech}{\mathrm{sech}}

\newcommand{\sign}{\mathop{\mathrm{sign}}}
\newcommand{\abs}[1]{\left\vert\,#1\,\right\vert}

\newcommand{\pd}[2]{\frac{\partial #1}{\partial\/ #2}}
\newcommand{\od}[2]{\frac{\mathrm{d} #1}{\mathrm{d}\/#2}}
\newcommand{\OD}[2]{\frac{\mathrm{D} #1}{\mathrm{D}\/#2}}
\newcommand{\norm}[1]{\left\vert\left\vert\,#1\,\right\vert\right\vert}

\newcommand{\half}{{\textstyle{1\over2}}}

\usepackage{acronym}
\acrodef{BVP}{Boundary Value Problem}
\acrodef{NSWE}{Nonlinear Shallow Water Equations}

\begin{document}

\title[\Title]{Numerical simulation of conservation laws with moving grid nodes: Application to tsunami wave modelling}

\author[G.~Khakimzyanov]{Gayaz Khakimzyanov}
\address{\textbf{G.~Khakimzyanov:} Institute of Computational Technologies of Siberian Branch of the Russian Academy of Sciences, Novosibirsk 630090, Russia}
\email{Khak@ict.nsc.ru}
\urladdr{https://www.researchgate.net/profile/Gayaz\_Khakimzyanov3/}

\author[D.~Dutykh]{Denys Dutykh$^*$}
\address{\textbf{D.~Dutykh:} Univ. Grenoble Alpes, Univ. Savoie Mont Blanc, CNRS, LAMA, 73000 Chamb\'ery, France and LAMA, UMR 5127 CNRS, Universit\'e Savoie Mont Blanc, Campus Scientifique, 73376 Le Bourget-du-Lac Cedex, France}
\email{Denys.Dutykh@univ-smb.fr}
\urladdr{http://www.denys-dutykh.com/}
\thanks{$^*$ Corresponding author}

\author[D.~Mitsotakis]{Dimitrios Mitsotakis}
\address{\textbf{D.~Mitsotakis:} Victoria University of Wellington, School of Mathematics, Statistics and Operations Research, PO Box 600, Wellington 6140, New Zealand}
\email{dmitsot@gmail.com}
\urladdr{http://dmitsot.googlepages.com/}

\author[N.~Yu.~Shokina]{Nina Yu. Shokina}
\address{\textbf{N.~Yu.~Shokina:} Department of Radiology, Medical Physics, Medical Center -- University of Freiburg, Faculty of Medicine, University of Freiburg, Freiburg, Germany}
\email{Nina.Shokina@gmail.com}
\urladdr{https://www.researchgate.net/profile/Nina\_Shokina/}


\begin{titlepage}
\thispagestyle{empty} 
\noindent
{\Large Gayaz \textsc{Khakimzyanov}}\\
{\it\textcolor{gray}{Institute of Computational Technologies, Novosibirsk, Russia}}
\\[0.02\textheight]
{\Large Denys \textsc{Dutykh}}\\
{\it\textcolor{gray}{CNRS, Universit\'e Savoie Mont Blanc, France}}
\\[0.02\textheight]
{\Large Dimitrios \textsc{Mitsotakis}}\\
{\it\textcolor{gray}{Victoria University of Wellington, New Zealand}}
\\[0.02\textheight]
{\Large Nina \textsc{Shokina}}\\
{\it\textcolor{gray}{Albert-Ludwigs-Universit\"at Freiburg, Freiburg im Breisgau, Germany}}\\[0.16\textheight]

\colorbox{Lightblue}{
  \parbox[t]{1.0\textwidth}{
    \centering\huge\sc
    \vspace*{0.7cm}
    
    \textcolor{bluepigment}{Numerical simulation of conservation laws with moving grid nodes: Application to tsunami \\ wave modelling}

    \vspace*{0.7cm}
  }
}

\vfill 

\raggedleft     
{\large \plogo} 
\end{titlepage}


\newpage
\maketitle
\thispagestyle{empty}


\begin{abstract}

In the present article we describe a few simple and efficient finite volume type schemes on moving grids in one spatial dimension combined with appropriate predictor--corrector method to achieve higher resolution. The underlying finite volume scheme is conservative and it is accurate up to the second order in space. The main novelty consists in the motion of the grid. This new dynamic aspect can be used to resolve better the areas with large solution gradients or any other special features. No interpolation procedure is employed, thus unnecessary solution smearing is avoided, and therefore, our method enjoys excellent conservation properties. The resulting grid is completely redistributed according the choice of the so-called \emph{monitor function}. Several more or less universal choices of the monitor function are provided. Finally, the performance of the proposed algorithm is illustrated on several examples stemming from the simple linear advection to the simulation of complex shallow water waves. The exact well-balanced property is proven. We believe that the techniques described in our paper can be beneficially used to model tsunami wave propagation and run-up.


\bigskip
\noindent \textbf{\keywordsname:} Conservation laws; finite volumes; conservative finite differences; moving grids; adaptivity; advection; shallow water equations; wave run-up \\

\smallskip
\noindent \textbf{MSC:} \subjclass[2010]{74S10 (primary), 74J15, 74J30 (secondary)}

\end{abstract}

\newpage
\tableofcontents
\thispagestyle{empty}


\newpage
\section{Introduction}

The main goal of the present manuscript consists in introducing to tsunami modelling community some state-of-the-art numerical methods for solving nonlinear hyperbolic equations, which allow to use the available degrees of freedom in the (nearly-) optimal way. Until the last $15$ years, tsunami wave modelling has been done in the framework of Nonlinear Shallow Water Equations (NSWE). The numerical simulations have been performed by practitioners essentially with several well-established community codes. To give a few examples, let us mention the code \textsc{TUNAMI} \cite{Imamura1996a} based on a conservative finite difference leap-frog scheme on staggered grids. This scheme approximates NSWE on real bathymetries. Another widely used code is \textsc{MOST}, which is based on the dimensional splitting method\footnote{This technique is also known as the alternating directions.} \cite{Titov1997, Titov1996}. We would like to mention also the software \textsc{COMCOT}, which is based on a modified leap-frog scheme discretizing NSWE in \textsc{Cartesian} and spherical coordinates \cite{Wang2011a}. Finally, in the Institute of Computational Technologies (SB RAS) an in-house software \textsc{MGC} based on the finite difference \textsc{MacCormack} scheme in spherical coordinates was developed as well \cite{Shokin2008}.

A big boost in tsunami research happened after Tsunami Boxing day \cite{Syno2006} due to the availability of real world data. After this \textsc{Sumatra} 2004 event, several studies questioned the importance of dispersive effects for trans-oceanic tsunami propagation \cite{DGK, Grilli2007, Murty2006}. Their conclusions were supported by comparisons of numerical results with field and DART buoy data. We can only mention that the importance of dispersive effects has been already highlighted in 1982 by \textsc{Mirchina} and \textsc{Pelinovsky} \cite{Mirchina1982}:
\begin{quote}
 \it [\,\dots\,] the considerations and estimates for actual tsunamis indicate that nonlinearity and dispersion can appreciable affect the tsunami wave propagation at large distances. [\,\dots\,]
\end{quote}

One of the first tsunami numerical models including spherical effects \emph{and} linear dispersive terms was \textsc{TUNAMI-N2}. A sensitivity study was performed in \cite{Tkalich2007} on the example of \textsc{Sumatra} 2004 event. The reported results were in favor of including dispersive effects. Another Weakly Nonlinear Weakly Dispersive (WNWD) model was presented in \cite{Lovholt2008, Lovholt2010}. It was basically a spherical counterpart of the classical \textsc{Peregrine} system \cite{Peregrine1967}. Let us mention also the study by \textsc{Glimsdal} \etal (2013) \cite{Glimsdal2013} where the importance of dispersion has been investigated numerically and the use of Fully Nonlinear Weakly Dispersive (FNWD) model was recommended. Such FNWD spherical equations (with the horizontal velocity defined on a certain surface inside the fluid bulk) have been derived in \cite{Kirby2013}, however, only WNWD numerical results were presented. In a recent work \cite{Khakimzyanov2016a} a systematic derivation of spherical rotating FNWD (and WNWD) models have been performed and the so-called \emph{base model} was derived, which incorporates many other models as particular cases. The reduction of spherical FNWD to WNWD models was illustrated and several known models were recovered in this way. Finally, the fully nonlinear numerical results have been presented in \cite{Khakimzyanov2016b}. In particular, the importance of sphericity, \textsc{Coriolis} and dispersion effects was thoroughly discussed. The numerical algorithm employed in \cite{Khakimzyanov2016b} is based on the splitting approach: an elliptic equation is solved to determine the non-hydrostatic pressure component and an evolution system of hyperbolic equations with source terms is solved to advance in time the velocity and water height variables. The algorithm was implemented as an explicit two-step predictor--corrector scheme. The volume of publication \cite{Khakimzyanov2016b} did not allow the authors to describe all details of the employed numerical method. Henceforth, if we adopt dispersive models for tsunami simulations\footnote{It seems now that future generation tsunami modelling codes will include at least some dispersive effects.}, there is a clear need in robust numerical methods to solve hyperbolic balance equations. The present manuscript should be considered as a tutorial-style article to introduce the moving grid algorithms to the developers of future generation tsunami propagation codes.

The numerical simulation of conservation laws is one of the most dynamic and central parts of modern numerical analysis. The systems of conservation laws appear in many fields ranging from traffic modeling \cite{Michalopoulos1993} to shallow water modelling \cite{Dutykh2010c} and compressible fluid mechanics \cite{Whitham1999}. The finite volume method was proposed in a pioneering work of S.~\textsc{Godunov} (1959) \cite{Godunov1959, Godunov1999} and developed later by Ph.~\textsc{Roe} \cite{Roe1981} and many other researchers (see \cite{Fuhrmann2014} for the current state of the art).

Nowadays the complexity of the problems which arise in practice is such that we have to think about the optimal usage of available computational resources. In this way the community came with the idea of developing numerical methods on adaptive grids in order to have higher resolution only where it is needed \cite{Babuska1995}. Historically, $hp$-adaptive methods were first developed for the Finite Element-type (FEM) discretizations \cite{Babuska1994} (including the newest discontinuous \textsc{Galerkin} discretizations as well \cite{Houston2001}). The main advantage of FEM is that this method is quite flexible in respect of local approximation spaces, thus allowing for relatively easy $p$-adaptivity. FEM was applied to hyperbolic problems under the guise of relaxation schemes \cite{Arvanitis2010}. Later $h$-adaptive methods have been adapted for finite volume methods as well \cite{Berger1984}. However, the widely used approach nowadays consists in performing local grid refinement (or coarsening) with locally nested grids \cite{Berger1989}. For instance, it has been successfully applied to tsunami propagation and run-up problems \cite{George2006a} as well as to more complex two-phase incompressible flows \cite{Popinet2003}.

Adaptive mesh refinement techniques applied to the spatial variable after the computation of an approximate solution of the model equations at each time step usually consist of the following steps:
\begin{itemize}
  \item Generation of a new spatial mesh according to prescribed adaptivity criteria;
  \item Reconstruction or interpolation of the numerical solution on the new mesh;
  \item Integration in time of the numerical solution on the new mesh.
\end{itemize}
The main disadvantages of this methodology is that the reconstruction of the numerical solution on the new mesh increases the complexity of the method and also reduces the accuracy of the method by adding more dissipation or dispersion to the solution. Moreover, the variable number of nodes represents some difficulties while implementing the algorithm on a computer. The underlying data structures have to be flexible enough to insert or to suppress some elements.

The approach we propose in this study is different and does not require interpolation or reconstruction of the solution to the new mesh (a conservative interpolation is employed in \cite{Tang2003}). Although the proposed method is closer to the redistribution method proposed in \cite{Tang2003, Arvanitis2006}, it is simpler and more elegant at the level of implementation details. Namely, it is based on the transformation of the model equations to some new equations whose solutions are independent of the new mesh. Namely, the method uses the same number of discretization points of the spatial mesh during the simulation and the adaptivity is achieved by moving the grid nodes to the places indicated by the so-called \emph{monitor function}. The adapted grid is obtained as a solution of an elliptic (or parabolic) problem, which ensures the smoothness of the obtained grid. The main idea behind the equidistribution principle is to distribute the nodes over the computational domain so that the measure of a cell times the value of the monitor function on it be approximatively constant \cite{Degtyarev1987a, Degtyarev1993}. The heart of the matter in the equidistribution method is the construction of a suitable monitor function. This method was described independently in the seminal paper \cite{Shokin1982}, and later in \cite{Huang1994}. Numerous subsequent developments were published in recent references \cite{Liseikin2010, Huang2011a, Gasparin2019a}. The review of earlier works on the adaptive grid generation can be found in \cite{Godunov1976, Liseikin1996}. The `Holy Grail' in choosing the monitor function is to achieve ideally the reduction of error in orders of magnitude when the discretization step $\Delta x\ \to\ 0\,$, as illustrated in \cite{Khakimzyanov2015b}. The main advantages of the proposed approach include:
\begin{itemize}
  \item Conservation in space;
  \item Second order accurate (in smooth regions);
  \item TVD property for scalar problems;
  \item Well-balanced character;
  \item Preservation of the discrete conservation law in cell centers and cell interfaces.
\end{itemize}
The present manuscript sheds also some light on the hyperobolic part of the splitting approach we used earlier to simulate fully nonlinear and weakly dispersive wave propagation \cite{Khakimzyanov2016}.

For simplicity, we focus on one-dimensional cases only, and even in this simple case some open problems are outlined throughout the manuscript. According to our knowledge the pioneers in 1D case were \cite{Sudobicher1968} (for nonlinear shallow water equations) and \cite{Alalykin1970} (for gas dynamics). The generation of 2D curvilinear adaptive grids is a relatively classical topic \cite{Serezhnikova1989}. A rigorous definition of a curvilinear grid was given in \cite{Shokin1985}. The generalization of the equidistribution method to two spatial dimensions was performed for the first time in \cite{Darmaev1987}. There are a few studies which report recent modern implementations on \emph{moving} redistributed two-dimensional grids (see \eg \cite{Cao1999, Azarenok2003, Tang2003, Beckett2002}) with generalized monitor functions formalized in \cite{Huang2011a}. Some approaches for the inter-comparison of generated grids are described in \cite{Prokopov1989}. We refer also to \cite{Budd2009a} for a recent review of moving grid techniques. The moving grid methods seem to be unavoidable in the modern numerical analysis. There are some classes of problems, which cannot be solved, in principle, with a fixed grid. To give an example, we can mention the class of problems involving moving boundaries such as free surface flows \cite{Bona2005, Knobloch2015, Khakimzyanov2018}.

The discretization of the model equations is based on explicit, fully discrete predictor-corrector Finite Volume schemes that are accurate up to the second order in space. Predictor--corrector schemes have been introduced by \textsc{MacCormack} in \cite{MacCormack1969} and can be considered as a natural generalization of splitting schemes in higher spatial dimensions.

We have to mention that non-uniform grids can be used also to preserve some \textsc{Lie} symmetry group of the equation at the discrete level (see \cite[Section~3.2.2]{Chhay2011} for a brief survey and the \textsc{Burgers}--\textsc{Hopf} equation example). Originally this approach was proposed by \textsc{Dorodnitsyn} \cite{Dorodnitsyn1994} and developed later in \cite{Huang2001}. In our study the mesh motion is directed solely by the equidistribution principle as in \cite{Tang2003} without taking into account the symmetry considerations for the moment.

The present article is organized as follows. In Section~\ref{sec:advect} we discuss in all details the implementation of the algorithm for the linear advection equation. Then, we generalize it to the nonlinear case in Section~\ref{sec:hopf} and systems of conservation laws in Section~\ref{sec:nswe}. Some open problems are outlined in Section~\ref{sec:open}. Finally, the main conclusions and perspectives of our study are given in Section~\ref{sec:concl}.


\section{Linear scalar equation}
\label{sec:advect}

In order to present the moving grid method we shall start with the simplest scalar linear advection equation in this section and, then, we shall increase gradually the complexity by adding first the nonlinearities in the following section and moving to the systems by the end of our manuscript.

We pay special attention to the linear case for the following reasons:
\begin{itemize}
  \item The main properties of the scheme are the most transparent in the linear case;
  \item The generalization to the nonlinear and vectorial cases will be easier once the linear case is fully understood. So, it will allow us to go faster in the subsequent sections;
  \item An exhaustive error analysis is possible in the linear (and presumably only in the linear) case.
\end{itemize}

In order to construct an efficient finite volume scheme on a moving mesh, we have to choose first a robust and an accurate scheme on a fixed grid, which will be generalized later to incorporate the motion of mesh points. Such a scheme retained for our study is described in the next Section.


\subsection{A predictor-corrector scheme on a fixed uniform mesh}
\label{sec:predcorr}

Consider the following linear advection equation with a constant propagation velocity
\begin{equation}\label{eq:adv}
  u_t\ +\ a\; u_x\ =\ 0, \qquad a\in \R.
\end{equation}
The subscripts in this study denote the partial derivatives, \ie $u_t = \pd{u}{t}$, $u_x = \pd{u}{x}$, \etc

We introduce a uniform discretization of the real space line $\R$ with nodes $x_j = j\,\Delta x$, where $\Delta x > 0$, $j\in\Z$ is the spatial discretization step, and for simplicity we do not pay attention to the boundary conditions here. The interval $\C_j = [x_j, x_{j+1}]$ will be referred as the cell $\C_j$. The time step is denoted by $\tau$ and the discrete solution is computed at $t^n = n\tau$, $n\in\Z^+$. Equation~\eqref{eq:adv} is discretized in space and time using an explicit fully discrete predictor-corrector scheme
\begin{equation}\label{eq:prcor2}
  \frac{u^*_{j+1/2} - u_{j+1/2}^n}{\tau^*_{j+1/2}}\ +\ a\; u_{x,j+1/2}^n\ =\ 0, \qquad 
  \frac{u_{j}^{n+1} - u_{j}^n}{\tau}\ +\ a\; \frac{u^*_{j+1/2} - u^*_{j-1/2}}{\Delta x}\ =\ 0,
\end{equation}
where $\tau$ is the time step. The intermediate quantities $u^*_{j\pm 1/2}$ are evaluated in the middle of cells at $x_{j\pm 1/2} = \displaystyle{\frac{x_j + x_{j \pm 1}}{2}} = x_j \pm \half \Delta x$ and at time instances $t = t^n + \tau^*_{j + 1/2}$. Moreover, the following notations have been introduced in \eqref{eq:prcor2}:
\begin{equation*}
  u_{j+1/2}^n\ =\ \frac{u_{j+1}^n + u_{j}^n}{2}, \quad
  u_{x,j+1/2}^n\ =\ \frac{u^n_{j+1} - u^n_{j}}{\Delta x}, \quad
  \tau^*_{j+1/2}\ =\ \frac{\tau}{2}\;\bigl(1\ +\ \theta^n_{j+1/2}\bigr).
\end{equation*}
We underline the fact that the parameter $\theta$ has not been fixed for the moment and it may vary from one cell and one time layer to another one.

The predictor-corrector scheme \eqref{eq:prcor2} can be recast as a one-step scheme by combining two equations together:
\begin{equation}\label{eq:can}
  \frac{u_{j}^{n+1} - u_{j}^{n}}{\tau}\ +\ a\;\frac{u_{j+1}^n - u_{j-1}^n}{2 \Delta x}\ -\ 
  \frac{\tau a^2}{2 \Delta x}\;\Bigl\{\bigl((1 + \theta)u_{x}\bigr)^n_{j+1/2}\ -\ \bigl((1 + \theta)u_{x}\bigr)^n_{j-1/2}\Bigr\}\ =\ 0.
\end{equation}
The last fully discrete scheme is a canonical form of all two steps explicit schemes for equation~\eqref{eq:adv}, since many well-known schemes can be obtained from \eqref{eq:can} by choosing carefully the parameter $\theta$. For example,
\begin{description}
  \item[Lax--Wendroff scheme \cite{Lax1960} ($\theta = 0$)]
  \begin{equation*}
    \frac{u_{j}^{n+1} - u_{j}^{n}}{\tau}\ +\ a\;\frac{u_{j+1}^n - u_{j-1}^n}{2 \Delta x}\ -\ \frac{\tau a^2}{2 \Delta x}\;\bigl\{u_{x,j+1/2}^n\ -\ u_{x,j-1/2}^n\bigr\}\ =\ 0.
  \end{equation*}
  \item[First order upwind \cite{Courant1952} ($\theta = \frac{1}{C} - 1 > 0$)]
  The number $C \coloneqq \abs{a}\;\dfrac{\tau}{\Delta x} < 1$ is the famous \textsc{Courant}--\textsc{Friedrichs}--\textsc{Lewy} number \cite{Courant1928}, which controls the stability of the upwind scheme
  \begin{equation*}
    \frac{u_{j}^{n+1} - u_{j}^n}{\tau}\ +\ \frac{a + \abs{a}}{2}\;u_{x,j-1/2}^n\ +\ 
    \frac{a - \abs{a}}{2}\;u_{x,j+1/2}^n\ =\ 0.
  \end{equation*}
  \item[Lax--Friedrichs scheme \cite{Breuss2004} ($\theta = \frac{1}{C^2} - 1 > 0$)]
  \begin{equation*}
    \frac{u^{n+1}_j\ -\ \frac{1}{2}(u_{j+1}^n + u_{j-1}^n)}{\tau}\ +\ a\;\frac{u_{j+1}^n - u_{j-1}^n}{2 \Delta x}\ =\ 0.
  \end{equation*}
\end{description}

However, the schemes listed above are not suitable for accurate simulations, since the first order upwind and \textsc{Lax}--\textsc{Friedrichs} schemes are too diffusive in practice and \textsc{Lax}--\textsc{Wendroff} scheme is not monotonicity preserving. A good solution consists in choosing the parameter $\theta$ adaptively. In a previous work \cite{Shokin2005} the following strategy was proposed:
\begin{equation*}
  \theta_{j+1/2}^n\ =\ \left\{
  \begin{array}{cl}
    0, & \abs{u_{x,j+1/2}^n}\ \leq\ \abs{u_{x,j+1/2-s}^n} \quad \mbox{and} \quad u_{x,j+1/2}^n \cdot u_{x,j+1/2-s}^n\ \geq\ 0, \\
    \theta_0 (1-\xi_{j+1/2}^n), & \abs{u_{x,j+1/2}^n}\ >\ \abs{u_{x,j+1/2-s}^n} \quad \mbox{and} \quad u_{x,j+1/2}^n \cdot u_{x,j+1/2-s}^n\ \geq\ 0, \\
    \theta_0, & u_{x,j+1/2}^n \cdot u_{x,j+1/2-s}^n\ <\ 0,
\end{array}\right.
\end{equation*}
with $\theta_0 \coloneqq \dfrac{1}{C} - 1 > 0$, $s \coloneqq \sign(a)$, $\xi_{j+1/2}^n \coloneqq \dfrac{u^n_{x,j+1/2-s}}{u^n_{x,j+1/2}}$. The first case corresponds to the classical \textsc{Lax}--\textsc{Wendroff} scheme, the third case is the classical first order upwind scheme and the second case is detailed in the following Remark~\ref{rem:1}. However, the proposed combination of schemes produces a robust TVD-type scheme, which ensures a very good trade-off between the solution accuracy and monotonicity properties. The scheme parameter $\theta$ plays somehow the r\^ole of the traditional \texttt{minmod} slope limiter in more conventional MUSCL-type finite volume methods \cite{LeVeque1992}. Its performance on fixed and mobile meshes will be illustrated in numerical sections below. The complete justification and validation of this scheme can be found in \cite{Shokin2005}.

\begin{remark}\label{rem:1}
We note that the canonical form \eqref{eq:can} contains not only three points schemes on a symmetric stencil. For instance, if we assume for simplicity $a > 0$, and set
\begin{equation*}
  \theta_{j+1/2}^n\ \coloneqq\ \Bigl(\frac{1}{C}\ -\ 1\Bigr)\Bigl(1\ -\ \frac{u_{x,j-1/2}^n}{u_{x,j+1/2}^n}\Bigr),
\end{equation*}
then we shall obtain the following second order upwind scheme on an asymmetric stencil:
\begin{equation*}
  \frac{u_{j}^{n+1} - u_{j}^{n}}{\tau}\ +\ a\;\frac{3u_{x,j-1/2}^n\ -\ u_{x,j-3/2}^n}{2}\ -\ \frac{\tau a^2}{2}\;\frac{u_{j}^{n}\ -\ 2u_{j-1}^n\ +\ u_{j-2}^n}{\Delta x^2}\ =\ 0.
\end{equation*}
\end{remark}

Of course, the numerical properties such as the accuracy, stability and monotonicity depend crucially on the choice of the parameter $\theta$ and have to be studied carefully in each particular situation.


\subsection{Predictor-corrector schemes on a moving grid}

In this section we shall restrict our attention to bounded domains only since almost all numerical simulations in 1D are made on finite intervals. To be more specific we assume that the computational domain is $\I\ =\ [\,0,\, \ell\,]\,$. Consider also the reference domain $Q\ =\ [\,0,\, 1\,]$ and a smooth bijective time-dependent mapping $x\,(\cdot,\, t):\ Q\,\times\,\R^{\,+}\ \mapsto\ \I\,$, which satisfies the following boundary conditions
\begin{equation}\label{eq:dir}
  x\,(0,\,t)\ =\ 0\,, \qquad x\,(1,\,t)\ =\ \ell\,.
\end{equation}
For our purely numerical purposes we do not even need to have the complete information about this mapping. Consider a uniform mesh of the reference domain $Q_{\,h}\ =\ \{q_{\,j}\ =\ j\,\Delta q\}_{\,j\,=\,0}^{\,N}\,$, with a constant step $\Delta q \coloneqq \dfrac{1}{N}\,$. Strictly speaking we need to know only the image of nodes $q_{\,j}$ under the map $x\,(q,\,t)\,$, since they constitute the nodes of the moving mesh:
\begin{equation*}
  x\,(q_{\,j},\, t^{\,n})\ =\ x_{\,j}^{\,n}\,.
\end{equation*}
The linear advection equation \eqref{eq:adv} can be rewritten on domain $Q$ if we introduce the composed function $v\,(q,\,t)\ \coloneqq\ (u\,\circ\, x)\,(q,\,t)\ \equiv\ u\,\bigl(x\,(q,\,t),\,t\bigr)\,$:
\begin{equation}\label{eq:trans}
  v_{\,t}\ +\ \frac{\ab}{J}\; v_q\ =\ 0\,,
\end{equation}
where $\ab\ \coloneqq\ a\ -\ x_{\,t}$ is the advection speed relative to the mesh nodes velocity and $J\ \coloneqq\ x_{\,q}\ >\ 0$ is the \textsc{Jacobian} of the transformation $x(q,t)$. Equation \eqref{eq:trans} does possess a conservative form as well
\begin{equation}\label{eq:cons}
  (Jv)_t\ +\ (\ab\, v)_q\ =\ 0.
\end{equation}
Now we can write the predictor-corrector scheme for the linear advection equation on a moving mesh:
\begin{align}\label{eq:twostages1}
  \frac{v^*_{j+1/2}\ -\ v^n_{j+1/2}}{\tau^*_{j+1/2}}\ +\ \Bigl(\frac{\ab}{J}\,v_q\Bigr)^n_{j+1/2}\ &=\ 0, \\
  \frac{(Jv)_{j}^{n+1}\ -\ (Jv)_{j}^n}{\tau}\ +\ \frac{\bigl(\ab^n\, v^*\bigr)_{j+1/2}\ -\ \bigl(\ab^n\, v^*\bigr)_{j-1/2}}{\Delta q}\ &=\ 0. \label{eq:twostages2}
\end{align}
Note that the first step is based on equation \eqref{eq:trans} and the second one on the conservative form \eqref{eq:cons}. Here we introduced some new notations:
\begin{equation*}
  x^n_{t,j+1/2}\ =\ \frac{x^n_{t,j}\ +\ x^n_{t,j+1}}{2}, \qquad 
  x^n_{t,j}\ =\ \frac{x_j^{n+1}\ -\ x_j^n}{\tau},
\end{equation*}
\begin{equation*}
  J_{j+1/2}^n\ \equiv\ x^n_{q,j+1/2}\ =\ \frac{x_{j+1}^n\ -\ x_{j}^n}{\Delta q}, \qquad 
  J_{j}^n\ \equiv\ x^n_{q,j}\ =\ \frac{J_{j+1/2}^n\ +\ J_{j-1/2}^n}{2}\ =\ \frac{x_{j+1}^n\ -\ x_{j-1}^n}{2 \Delta q}.
\end{equation*}
As before, one can recast the scheme \eqref{eq:twostages1}, \eqref{eq:twostages2} in an equivalent one-stage form:
\begin{equation*}
  \frac{(Jv)_{j}^{n+1}-(Jv)_{j}^n}{\tau}\ +\ \frac{\left({\bar a} v \right)^n_{j+1/2} - \left({\bar a} v \right)^n_{j-1/2}}{\Delta q}\ -\ \frac{\tau}{2\Delta q} \biggl\{\Bigl((1+\theta)\frac{\ab^2}{J}v_q\Bigr)^n_{j+1/2} - \Bigl((1+\theta) \frac{\ab^2}{J}v_q\Bigr)^n_{j-1/2}\biggr\}\ =\ 0.
\end{equation*}
A stability study using the method of frozen coefficients \cite{Rozhdestvenskiy1978} leads to the following restrictions on the scheme parameters
\begin{equation}\label{eq:stab}
  \theta_{j+1/2}^n \geq 0, \qquad \max\limits_{j\in \{0,\ldots,N-1\}} \bigl(\sqrt{1+\theta}\; C\bigr)_{j+1/2}^n \leq 1,
\end{equation}
where $C_{j+1/2}^n$ is the local CFL number \cite{Courant1928} defined as
\begin{equation*}
  C_{j+1/2}^n\ =\ \frac{\tau}{\Delta q}\; \biggl(\frac{\abs{\ab}}{J}\biggr)^n_{j+1/2}.
\end{equation*}
In order to ensure the monotonicity of the numerical solution, it was shown in \cite{Shokin2005} that it is sufficient to use the following choice of the scheme parameter $\theta$:
\begin{equation}\label{eq:21}
  \theta_{j+1/2}^n\ =\ \left\{
  \begin{array}{cl} 
  0, & \abs{\gt_{j+1/2}^n}\ \leq\ \abs{\gt_{j+1/2-s}^n} \quad \mbox{and} \quad \gt_{j+1/2}^n \cdot \gt_{j+1/2-s}^n\ \geq\ 0, \\
  \bigl(\theta_{0}(1 - \xi)\bigr)_{j+1/2}^n, & \abs{\gt_{j+1/2}^n}\ >\ \abs{\gt_{j+1/2-s}^n} \quad \mbox{and} \quad \gt_{j+1/2}^n \cdot \gt_{j+1/2-s}^n\ \geq\ 0, \\
  \theta_{0,j+1/2}^n, & \gt_{j+1/2}^n \cdot \gt_{j+1/2-s}^n\ <\ 0,
\end{array}\right.
\end{equation}
where $s = \sign(\ab_{j+1/2}^n)$ and
\begin{equation*}
  \theta_{0,j+1/2}^n\ =\ \frac{1}{C_{j+1/2}^n}\ -\ 1, \qquad
  \gt_{j+1/2}^n\ =\ \bigl(\abs{\ab}\;(1 - C)\; v_q\big)_{j+1/2}^n, \qquad
  \xi_{j+1/2}^n\ =\ \frac{\gt_{j + 1/2 - s}^n}{\gt_{j + 1/2}^n}.
\end{equation*}
Using \eqref{eq:21}, it can be shown that stability requirements \eqref{eq:stab} are fulfilled under the following single restriction on local CFL numbers:
\begin{equation*}
  C_{j+1/2}^n\ \leq\ 1.
\end{equation*}
Finally, the global CFL number $C$ is defined as
\begin{equation*}
  C\ \coloneqq\ \max_{j}\{C_{j+1/2}^n\}.
\end{equation*}
It is this value, which will be reported in Tables below.

There is an additional relation, which is called the \emph{geometric conservation law}. In 1D case it takes a very simple form
\begin{equation*}
  J_{\,t}\ -\ (x_{\,t})_{\,q}\ =\ 0\,,
\end{equation*}
since it translates a trivial fact that $x_{\,q\,t}\ \equiv\ x_{\,t\,q}\,$. However, it was demonstrated in \cite{Thomas1979} that it is absolutely crucial to respect the geometric conservation law at the discrete level as well in order to have a consistent fully discrete scheme which preserves all constant solutions exactly. It can be shown that the scheme proposed above satisfies the discrete counterpart of the geometric conservation law in integer and half-integer nodes:
\begin{align*}
  \frac{J^{n+1}_{j}\ -\ J^{n}_{j}}{\tau}\ -\ \frac{x^n_{t,j+1/2}\ -\ x^n_{t,j-1/2}}{\Delta q}\ &=\ 0, \\
  \frac{J^{n+1}_{j+1/2}\ -\ J^{n}_{j+1/2}}{\tau}\ -\ \frac{x^n_{t,j+1}\ -\ x^n_{t,j}}{\Delta q}\ &=\ 0.
\end{align*}


\subsection{Construction and motion of the grid}

In the previous sections we explained in some detail the predictor-corrector schemes on fixed and moving meshes. We saw that the mesh motion is parametrized by a bijective mapping $x(q,t): Q\times\R^+\mapsto \I$. Now we explain how this map is effectively constructed, since so far it was quite general. Moreover, the predictor-corrector schemes presented above are valid for other constructions of the mesh motion. 

The first step consists in choosing a positive valued function $w\,(x,\,t)\,:\ \I\,\times\,\R^{\,+}\ \mapsto\ \R^{\,+}\,$, the so-called \emph{monitor function}, which controls the distribution of the nodes. Without entering into the details for the moment, we suggest as one choice for a monitor function the following expression:
\begin{equation}\label{eq:def}
  w\,(x,\,t)\ =\ 1\ +\ \alpha\;\abs{u_{\,x}\,(x,t)}\,,
\end{equation}
where $\alpha\ >\ 0$ is a real parameter. The monitor function \eqref{eq:def} is a particular case of the class of monitor functions considered by \textsc{Beckett}, \textsc{Mackenzie} \etal \cite{Beckett2002}. Other choices of the monitor function are discussed in \cite{Budd2009, Stockie2001, Cao1999}. For instance, to give an example, another popular choice of the monitor function is 
\begin{equation*}
  w\,(x,\,t)\ =\ \sqrt{1\ +\ \alpha\;\abs{u_{\,x}\,(x,\,t)}^{\,2}}\,.
\end{equation*}
In general, the choice of the monitor function and its parameters might be crucial for the accuracy of the scheme. The adaptive grid construction for simple linear singularly perturbed \textsc{Sturm}--\textsc{Liouville} equation is deeply analysed in \cite{Khakimzyanov2015b}. This question seems to be rather open for unsteady nonlinear problems. The current state-of-the-art in this field is described in \cite{Huang2011a}.

It is noted that the monitor function varies in both space and time and takes large (positive) values in areas where the solution has important gradients. The monitor function $w(x,t)$ specified in \eqref{eq:def} can be readily discretized
\begin{equation*}
  w^{\,n}_{\,j\,+\,1/2}\ =\ 1\ +\ \alpha\;\frac{\abs{u_{\,j\,+\,1}^{\,n}\ -\ u_{\,j}^{\,n}}}{x_{\,j\,+\,1}^{\,n}\ -\ x_{\,j}^{\,n}}\,.
\end{equation*}

The proposed algorithm is slightly different on the very first step when the initial grid is generated. It is important to obtain a high-quality initial mesh, since the numerical error committed in the beginning cannot be corrected afterwards. So, the first step is explained in detail below and the algorithm for unsteady computations is presented in the following section.


\subsubsection{Initial grid generation}

It is of capital importance to produce a grid adapted to the initial condition as well, since no error committed initially can be corrected during dynamical simulations. At $t\ =\ 0$ we compute the monitor function $w\,(x,\,0)$ of the initial condition $u\,(x,\,0)$ and the mapping $x\,(q,\,0)$ is determined as the solution to the following \emph{nonlinear} elliptic problem, which degenerates to a simple second order ordinary differential equation
\begin{equation}\label{eq:steady}
  \bigl(\,w\,(x,\,0)\;x_{\,q}\,\bigr)_{\,q}\ =\ 0\,,
\end{equation}
supplemented with \textsc{Dirichlet} boundary conditions \eqref{eq:dir}. This elliptic problem can be readily discretized to produce the following finite difference approximation:
\begin{equation*}
  \frac{1}{\Delta q}\;\biggl\{w^0_{j+1/2}\;\frac{x^0_{j+1}\ -\ x^0_j}{\Delta q}\ -\ w^0_{j-1/2}\;\frac{x^0_j\ -\ x^0_{j-1}}{\Delta q}\biggr\}\ =\ 0, \quad j\ =\ 1, \ldots, N - 1,
\end{equation*}
with boundary conditions $x_0^0\, =\, 0$, $x_N^0\, =\, \ell\,$. The last nonlinear system of equations is solved iteratively (usually with some fixed point-type iterations). The iterations are initialized with a uniform grid as the first guess. One can notice also that its solution satisfies the equidistribution principle:
\begin{equation}\label{eq:equi}
  w^0_{j+1/2}\;\frac{x_{j+1}^0\ -\ x_j^0}{\Delta q}\ =\ C_{\Delta q}, \qquad j\ =\ 0,\ldots, N - 1,
\end{equation}
where the constant $C_h$ can be computed as
\begin{equation*}
  C_{\Delta q}\ =\ \sum_{j=0}^{N-1} (\Delta x\cdot w)_{j+1/2}^0, \qquad
  \Delta x^0_{j+1/2}\ \coloneqq\ x_{j+1}^0\ -\ x_j^0.
\end{equation*}
Therefore, $C_{\Delta q}$ is the discrete version of the monitor function integral. The equidistribution principle gives meaning to the monitor function: in the areas where $w^0_{j+1/2}$ takes large values, the spacing between two neighboring nodes $x_{j+1}^0$ and $x_j^0$ has to be inversely proportionally small.


\subsubsection{Unsteady computations}

Assume that the grid $\{x_j^n\}_{j=0}^N$ at $t\, =\, t^n$ is known and let us evolve it to the next time layer $t^{n+1}$. It is done by solving the following \emph{linear} fully discrete problem:
\begin{equation}\label{eq:diff}
  \frac{1}{\Delta q}\;\biggl\{w^n_{j+1/2}\;\frac{x^{n+1}_{j+1}\ -\ x^{n+1}_j}{\Delta q}\ -\ w^n_{j-1/2}\;\frac{x^{n+1}_j\ -\ x^{n+1}_{j-1}}{\Delta q}\biggr\}\ =\ \beta\;\frac{x^{n+1}_j - x^{n}_j}{\tau}, \quad j\ =\ 1, \ldots, N - 1,
\end{equation}
with the same boundary conditions $x_0^{n+1} = 0$, $x_N^{n+1} = \ell$ as above. Here the real positive parameter $\beta > 0$ plays the r\^ole of the diffusion coefficient and it controls the smoothness of nodes trajectories. The \emph{optimal} choice of the parameter $\beta$ is currently an open problem and it has to be done empirically in case by case basis. It can be easily seen that the scheme \eqref{eq:diff} is an implicit discretization of this nonlinear parabolic equation:
\begin{equation}\label{eq:parab}
  \bigl(\,w\,(x,\,t)\,x_{\,q}\,\bigr)_{\,q}\ =\ \beta\,x_{\,t}\,.
\end{equation}
However, it has to be stressed out that the fully discrete problem is linear since the monitor function $w^n_{j \pm 1/2}$ is taken on the previous time layer. Thus, the overall complexity of the grid motion algorithm at every time step is equal to that of the tridiagonal matrix algorithm \cite{Fedorenko1994}.


\subsubsection{A lyrical digression}

It is probably not completely clear to the reader why one has to use a different procedure initially, especially since the solution at every time step might be considered as initial condition to the next time level. This question would be completely legitimate.

The moving grid technique based on the solution of a (nonlinear) parabolic Equation~\eqref{eq:parab} was specifically designed for the solution of unsteady problems. It was proposed for the first time in \cite{Shokin1982} to simulate gas dynamics problems and it allows to ensure the required \emph{smoothness} of grid points trajectories (see numerical examples below). The position of grid points on two consecutive time steps are related by a finite-difference analogue of Equation~\eqref{eq:parab}. This modification of the equidistribution principle allows to avoid any rapid changes in the location of grid nodes. Basically, it was understood in the original work \cite{Shokin1982} and confirmed later in subsequent developments in late (19)80's and early (19)90's \cite{Dorfi1987, Darin1988, Degtyarev1993, Zegeling1992, Zegeling1992a}. Subsequent developments have been done in \cite{Huang1994, Huang2001, Huang2011a}. Some recommendations against using the steady Equation~\eqref{eq:steady} to determine the position of grid nodes on each time step were formulated in \cite{Daripa1992, Zegeling1992a, Coyle1986}. Some examples of grid nodes oscillations and other instabilities were shown there as well. In Appendix~\ref{app:ex} we provide a few \emph{analytical} examples to illustrate eventual problems which may arise if one uses the equidistribution principle \eqref{eq:steady} on every time step. Indeed, we demonstrate in Appendix~\ref{app:ex} that even infinitesimal perturbations of the monitor function may yield finite perturbations in the solution $x\,(q,\,t)\,$.


\subsubsection{Smoothing step}

In some problems where the solution contains a lot of oscillations (\ie many local extrema), the monitor function $w(x,\,t)$ has to be filtered out to ensure the smoothness of the mesh motion. Moreover, a smoothing step allows to enlarge the set of acceptable values of the parameter $\alpha$ in the definition \eqref{eq:def} of the monitor function. Our numerical experiments show that the following implicit procedure (inspired by implicit discretizations of the heat equation) produces very robust results:
\begin{equation}\label{eq:sys}
  \wb_{j+1/2}\ =\ w_{j+1/2}\ -\ \sigma\wb_{j+1/2}\ +\ \frac{\sigma}{2}\;(\wb_{j-1/2}\ +\ \wb_{j+3/2})\,, \qquad j\ =\ 1,\ldots,N-2\,,
\end{equation}
where $\sigma > 0$ is an ad-hoc smoothing parameter to be set empirically. In our study we do smoothing similarly to \textsc{Asselin}'s filter proposed earlier for time integration problems \cite{Asselin1972}. System \eqref{eq:sys} is completed with boundary conditions
\begin{equation*}
  \wb_{1/2}\ =\ w_{1/2}\,, \qquad \wb_{N-1/2}\ =\ w_{N-1/2}\,.
\end{equation*}
The smoothed discrete monitor function $\{\wb_j\}_{j=0}^N$ is obtained by solving the linear system \eqref{eq:sys} with any favorite method. In the present study we used the simplest tridiagonal matrix algorithm \cite{Fedorenko1994}. In contrast, the authors of \cite{Tang2003} used an explicit smoothing method. We stress out that one should apply a smoothing operator to the monitor function and not to the numerical solution.

\subsubsection{On the choice of the monitor function}

In general, when one has a problem (\ie a system of PDEs) to solve, there is a question of the best possible monitor function choice for the grid adaptation. One has to choose also the pertinent solution component as an input for the monitor function. It is extremely difficult to give some general recommendations which would work in all cases\footnote{It is probably as difficult as the general theory of PDEs.}. However, in choosing the monitor function one has to take into account the following aspects:
\begin{description}
  \item[Physical considerations] any a priori knowledge on the solution can be used to underline important physical aspects (\eg the presence of shock waves or the presence of multiple wave crests, \etc);
  \item[Geometrical considerations] computational domain geometry (\eg the presence of angles or fractal shorelines in wave run-up problems) and topology (\ie the presence of islands or other obstacles);
  \item[Mathematical considerations] nature of equations (differential, partial differential, integro-differential, singularly perturbed, \etc), but also one has to pick up the most pertinent variable to compute the monitor function on it (the free surface elevation or the vorticity field in \textsc{Navier--Stokes} simulations);
  \item[Numerical considerations] particularities of the employed discretization scheme may also hint the choice of the monitor function.
\end{description}
Once the monitor function is chosen, one has to determine also (nearly) optimal values of parameters which dependent not only on the problem, but also on the solution. Finally, the most important guideline nowadays is the \emph{experience} in solving practical problems using moving grid techniques. It cannot be replaced by any guidelines. We refer also to the review \cite{Budd2009a}.

\subsection{Numerical results}

\begin{figure}
  \centering
  \subfigure[]{\includegraphics[width=0.49\textwidth]{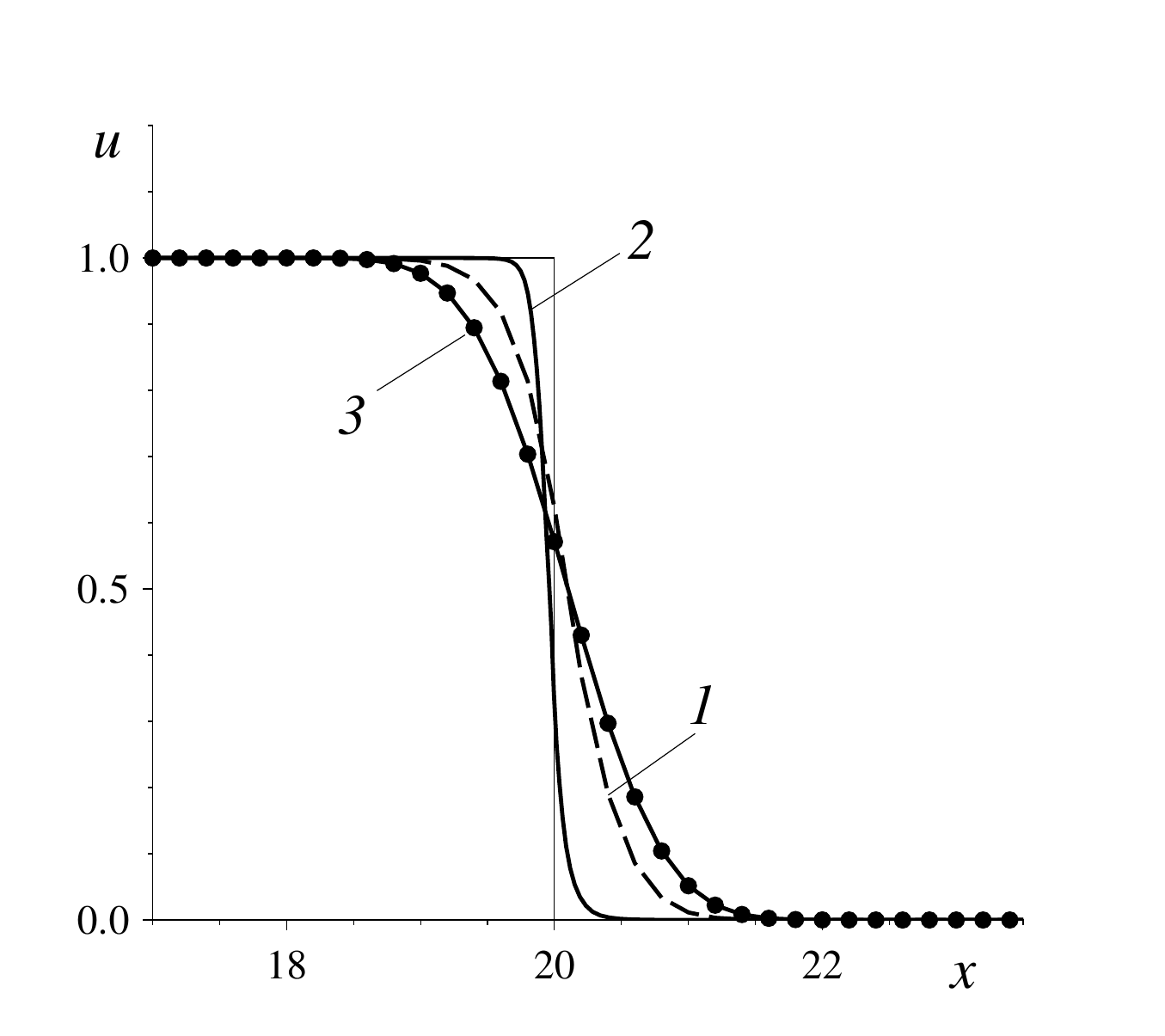}}
  \subfigure[]{\includegraphics[width=0.49\textwidth]{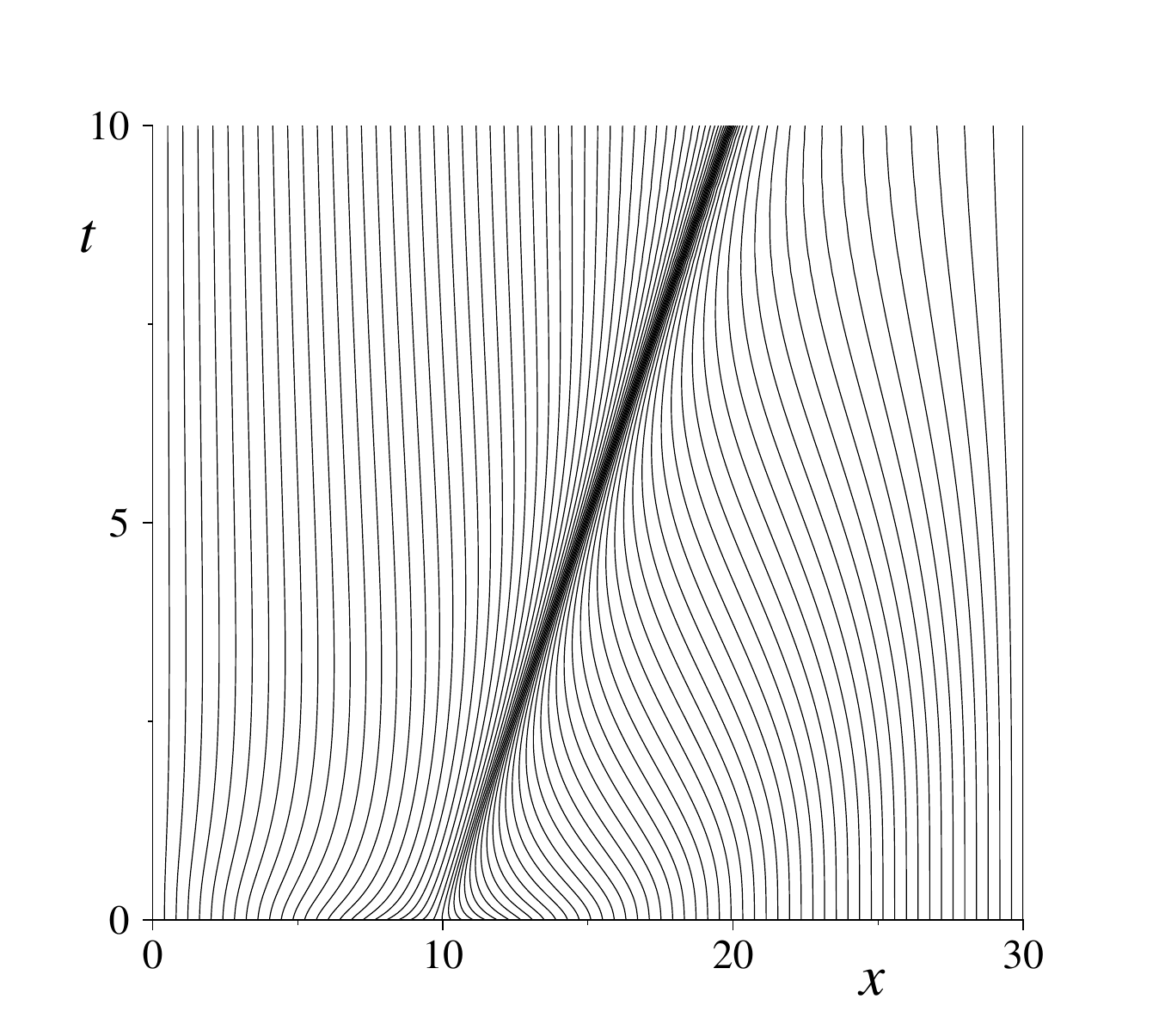}}
  \caption{\small\em Discontinuous profile propagation under the linear advection equation \eqref{eq:adv}: (a) comparison of different profiles (thin solid line is the exact solution, (1) is the numerical solution on the uniform mesh, (2) is the solution on the moving mesh and (3) is the first order upwind scheme on the uniform mesh); (b): trajectories of grid nodes in space-time. Note that the panel (a) shows only a zoom of the computational domain.}
  \label{fig:lshock}
\end{figure}

In order to illustrate the performance of this algorithm we choose a notorious test case, which consists in simulating the propagation of a discontinuous profile: 
\begin{equation*}
  u(x,0)\ =\ \left\{\begin{array}{ll}
    1, & x\ \leq\ x_*, \\
    0, & x\ >\ x_*,
  \end{array}\right.
  \qquad x_*\in ]\,0, \ell\,[\,.
\end{equation*}
The monitor function employed in this computation is defined in \eqref{eq:def}. The numerical parameters used in this simulations are given in Table~\ref{tab:params1}. The simulation results at time $t\ =\ T_{\,f}$ is shown in Figure~\ref{fig:lshock}. From the comparison of profiles obtained on fixed and moving meshes, one can see that the moving mesh reduces significantly the ``\emph{width}'' of the numerical shock wave profile.

\begin{table}
  \centering
  \begin{tabular}{lc}
  \hline\hline
  \textit{Parameter} & \textit{Value} \\
  \hline
  Computational domain length, $\ell$ & $30.0$ \\
  Number of nodes, $N$ & $150$ \\
  CFL number, $C$ & $0.8$ \\
  Final simulation time, $T_f$ & $10.0$ \\
  Propagation speed, $a$ & $1.0$ \\
  Initial position of the discontinuity, $x_*$ & $10.0$ \\
  Monitor function parameter, $\alpha$ & $10.0$ \\
  Grid diffusion parameter, $\beta$ & $150.0$ \\
  Smoothing parameter, $\sigma$ & $100.0$ \\
  \hline\hline
  \end{tabular}
  \bigskip
  \caption{\small\em Numerical parameters used to simulate the discontinuous profile propagation under the linear advection equation \eqref{eq:adv}.}
  \label{tab:params1}
\end{table}

The second test case is a smooth bell-shaped function advected by the dynamics of \eqref{eq:adv}. This test is used to measure the ability of the scheme to preserve the extrema (in other words we assess the strength of the numerical diffusion on a given mesh). The exact solution is given by 
\begin{equation*}
  u(x,t)\ =\ \ue^{-25(x\ -\ x_0\ -\ a t)^2}.
\end{equation*}
Thus, the initial condition is simply given by $u(x,0) = \ue^{-25(x - x_0)^2}$. Another particularity of the present test case is that we use the monitor function specifically designed with emphasis on the fine resolution of the local extrema (those with high positive or high negative values):
\begin{equation}\label{eq:defd2}
  w(x,t)\ =\ 1\ +\ \alpha\abs{u(x,t)}.
\end{equation}
It is noted that there is no differentiation operator in the definition of the last monitor function compared to the one defined in (\ref{eq:def}).

\begin{figure}
  \centering
  \subfigure[]{\includegraphics[width=0.49\textwidth]{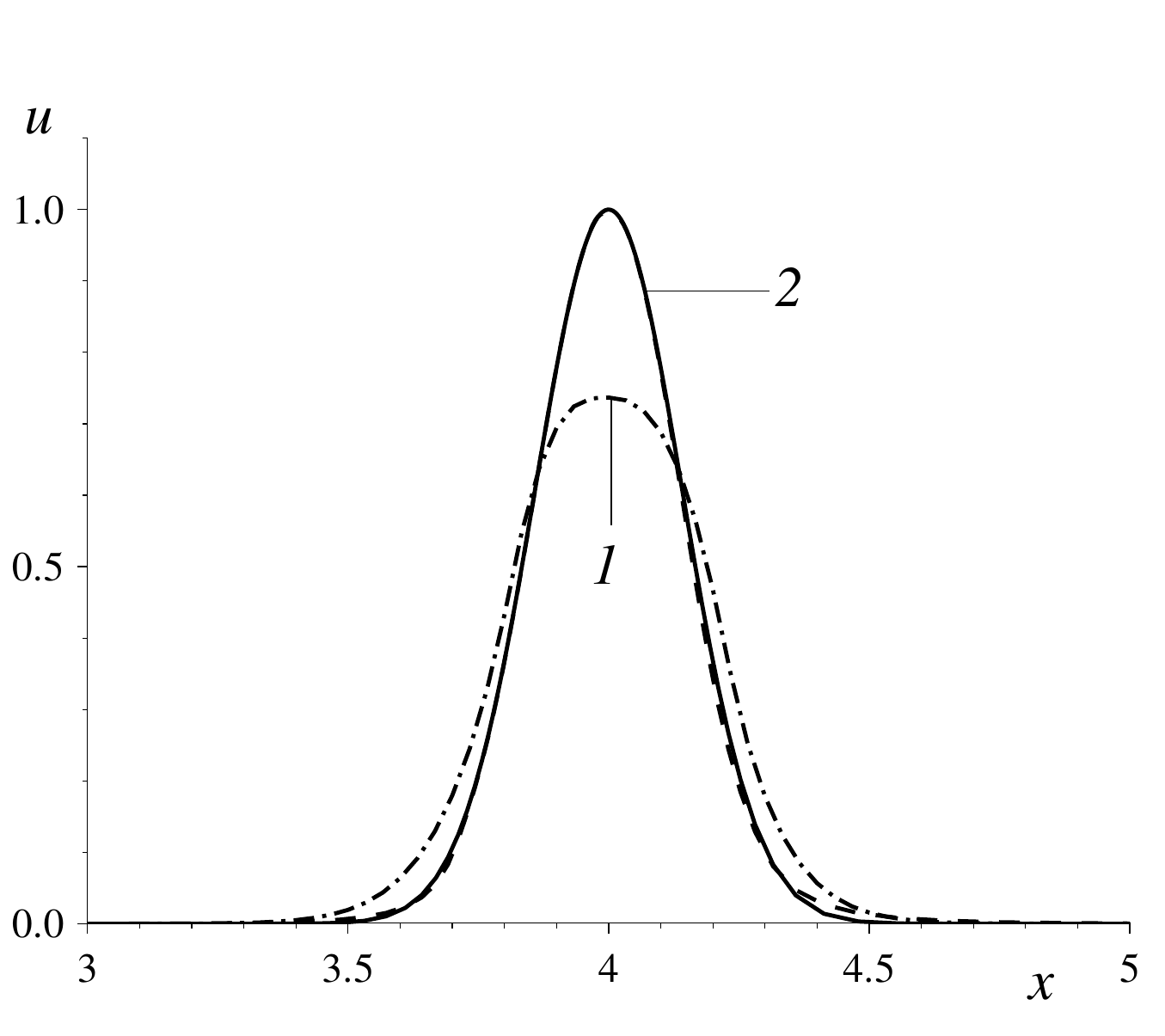}}
  \subfigure[]{\includegraphics[width=0.49\textwidth]{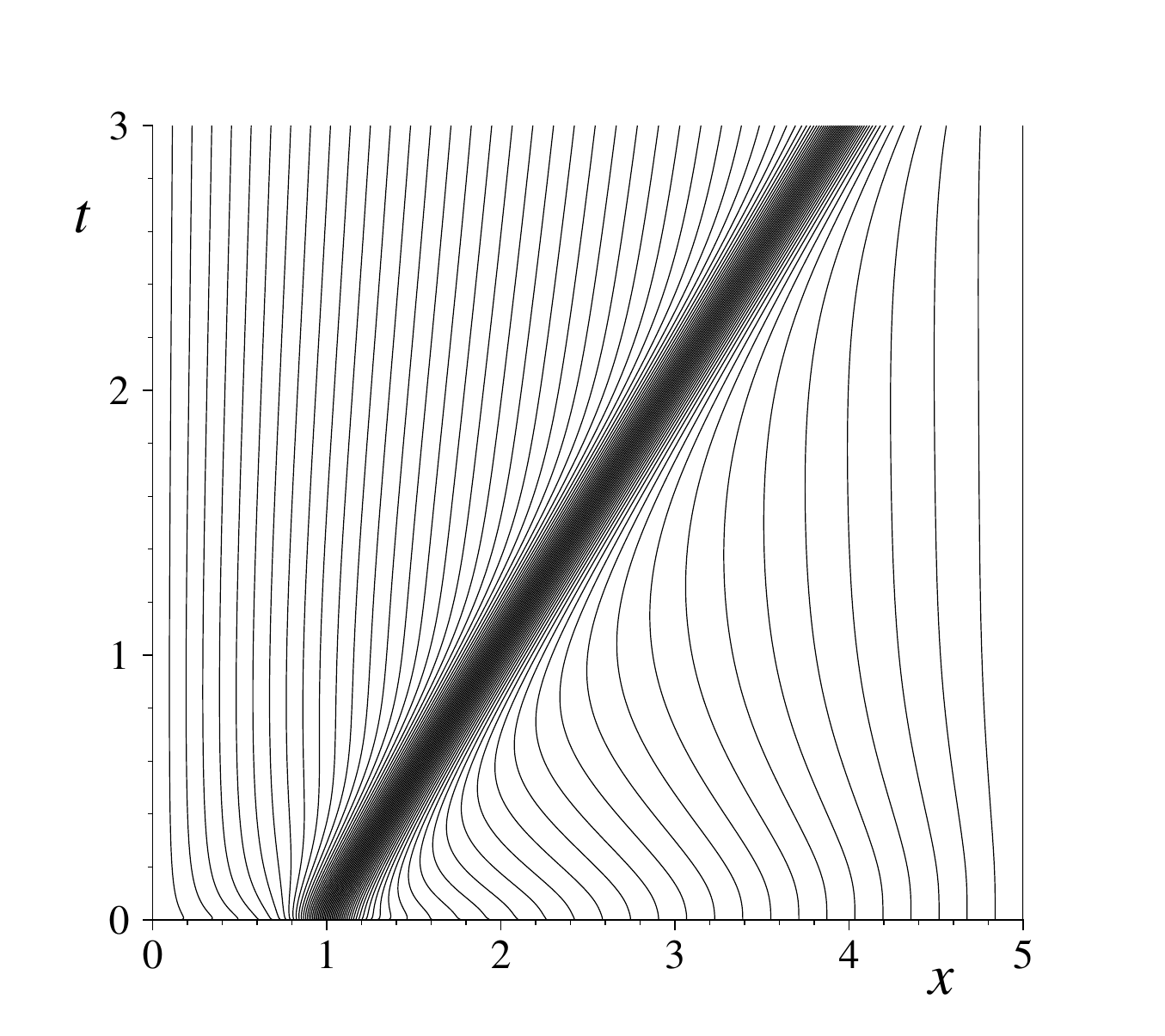}}
  \caption{\small\em Smooth profile propagation under the linear advection equation \eqref{eq:adv}: (a) comparison of different profiles (thin solid line is the exact solution, (1) is the numerical solution on the uniform mesh, (2) is the solution on the moving mesh); (b): trajectories of grid nodes in space-time. Note that the panel (a) shows only a zoom of the computational domain.}
  \label{fig:lexp}
\end{figure}

The parameters used in the numerical simulation are provided in Table~\ref{tab:params2}. The numerical results at $t = T_f$ are shown in Figure~\ref{fig:lexp}. For the sake of comparison we show also the numerical solution obtained with the same scheme, but on a uniform mesh. The use of a moving grid improves dramatically the preservation of the local maximum. One cannot even distinguish the exact solution from the numerical one up to the graphic resolution.

\begin{table}
  \centering
  \begin{tabular}{lc}
  \hline\hline
  \textit{Parameter} & \textit{Value} \\
  \hline
  Computational domain length, $\ell$ & $5.0$ \\
  Number of nodes, $N$ & $150$ \\
  CFL number, $C$ & $0.8$ \\
  Final simulation time, $T_f$ & $3.0$ \\
  Initial position of the bell, $x_0$ & $1.0$ \\
  Propagation speed, $a$ & $1.0$ \\
  Monitor function parameter, $\alpha$ & $20.0$ \\
  Grid diffusion parameter, $\beta$ & $20.0$ \\
  Smoothing parameter, $\sigma$ & $10.0$ \\
  \hline\hline
  \end{tabular}
  \bigskip
  \caption{\small\em Numerical parameters used to simulate the bell-shaped profile propagation under the linear advection equation \eqref{eq:adv}.}
  \label{tab:params2}
\end{table}


\section{Nonlinear scalar equation}
\label{sec:hopf}

In this section we consider a general hyperbolic scalar nonlinear conservation law of the form:
\begin{equation*}
  u_t\ +\ [f(u)]_x\ =\ 0.
\end{equation*}
As a chrestomathic example of such equations we can mention the celebrated \textsc{Burgers}--\textsc{Hopf} equation \cite{Burgers1948}. As above, we make a change of variables to the fixed reference domain $Q$:
\begin{equation*}
  v_t\ +\ \frac{1}{J}\;\bigl[f_q\ -\ x_t\,v_q\bigr]\ =\ 0,
\end{equation*}
or in the conservative form:
\begin{equation*}
  (Jv)_t\ +\ \bigl[f\ -\ x_t\,v\bigr]_q\ =\ 0,
\end{equation*}
where $v \coloneqq u\circ x$ and $J \coloneqq x_q$ as before.

For nonlinear equations we need an extra transport equation for the flux function:
\begin{equation*}
  f_t\ +\ a(u)\;f_x\ =\ 0,
\end{equation*}
where $a(u) = f_u(u)$. The last equation transforms to the following non-conservative form on the fixed domain $Q$:
\begin{equation*}
  f_t\ +\ \frac{a(v)}{J}\;\bigl[f_q\ -\ x_t\,v_q\bigr]\ =\ 0.
\end{equation*}
Now we can write down the predictor-corrector scheme directly on a moving grid (notice that in the nonlinear case three stages are required):
\begin{align*}
  \frac{v^*_{j+1/2}\ -\ v_{j+1/2}^n}{\tau^*_{j+1/2}}\ +\ \biggl(\frac{f_q\ -\ x_t\,v_q}{J}\biggr)_{j+1/2}^n\ &=\ 0, \\
  \frac{f^*_{j+1/2}\ -\ f_{j+1/2}^n}{\tau^*_{j+1/2}}\ +\ \biggl(\frac{a(f_q\ -\ x_t\,v_q)}{J}\biggr)_{j+1/2}^n\ &=\ 0, \\
  \frac{(Jv)_{j}^{n+1}\ -\ (Jv)_{j}^n}{\tau}\ +\ \frac{(f^*\ -\ x^n_t\,v^*)_{j+1/2}\ -\ (f^*\ -\ x^n_t\,v^*)_{j-1/2}}{\Delta q}\ &=\ 0,
\end{align*}
where $v_{j+1/2}^n$, $f_{j+1/2}^n$ were defined above. The last scheme has the second order approximation in space $\O(\Delta q^2)$ for smooth solutions provided that $\theta = \O(\Delta q)$. The only point which remains to be specified is the computation of the numerical \textsc{Jacobian} $a_{j + 1/2}^n$. In order to ensure good properties of the scheme, one has to ensure the compatibility condition $f_q = a\,(v)\,v_{\,q}$ at the discrete level $f^n_{q,j+1/2} = (a\,v_q)^n_{j+1/2}$ for finite differences analogue of the derivatives:
\begin{equation*}
  f^n_{q,j+1/2}\ =\ \frac{f^n_{j+1}\ -\ f^n_{j}}{\Delta q}, \qquad
  v^n_{q,j+1/2}\ =\ \frac{v^n_{j+1}\ -\ v^n_{j}}{\Delta q}.
\end{equation*}
This requirement can be satisfied with the following choice \cite{Harten1983}:
\begin{equation*}
  a^n_{j+1/2}\ =\ \left\{
  \begin{array}{cl}
    \frac{f_{j+1}^n\ -\ f_{j}^n}{v_{j+1}^n\ -\ v_{j}^n}, & v_{j+1}^n\ \neq\ v_{j}^n, \\
    a(v_{j}^n), & v_{j+1}^n\ =\ v_{j}^n.
\end{array}\right.
\end{equation*}

By noticing that the corrector step contains only the combination $f^*-x_t^n\,u^*$, we can write a compact two-stage form of the predictor-corrector scheme:
\begin{align*}
  \Bigl(\frac{\fh\ -\ f^n}{\tau^*}\Bigr)_{j+1/2}\ +\ \Bigl(\frac{\ab^2}{J}\;v_q\Bigr)^n_{j+1/2}\ &=\ 0, \\
  \frac{(Jv)_{j}^{n+1}\ -\ (Jv)_{j}^n}{\tau}\ +\ \frac{\bigl(\fh\ -\ (x_{t}\,v)^n\bigr)_{j+1/2}\ -\ \bigl(\fh\ -\ (x_{t}\,v)^n\bigr)_{j-1/2}}{\Delta q}\ &=\ 0,
\end{align*}
where $\ab(v,q,t) = a(v) - x_t$. The parameter $\theta_{j+1/2}^n$ is computed according to \eqref{eq:21} which ensures the TVD property of the scheme. As a result we obtain a generalization of the \textsc{Harten} scheme \cite{Harten1983} to the case of moving meshes.

\subsection{Numerical results}

\begin{figure}
  \centering
  \subfigure[]{\includegraphics[width=0.49\textwidth]{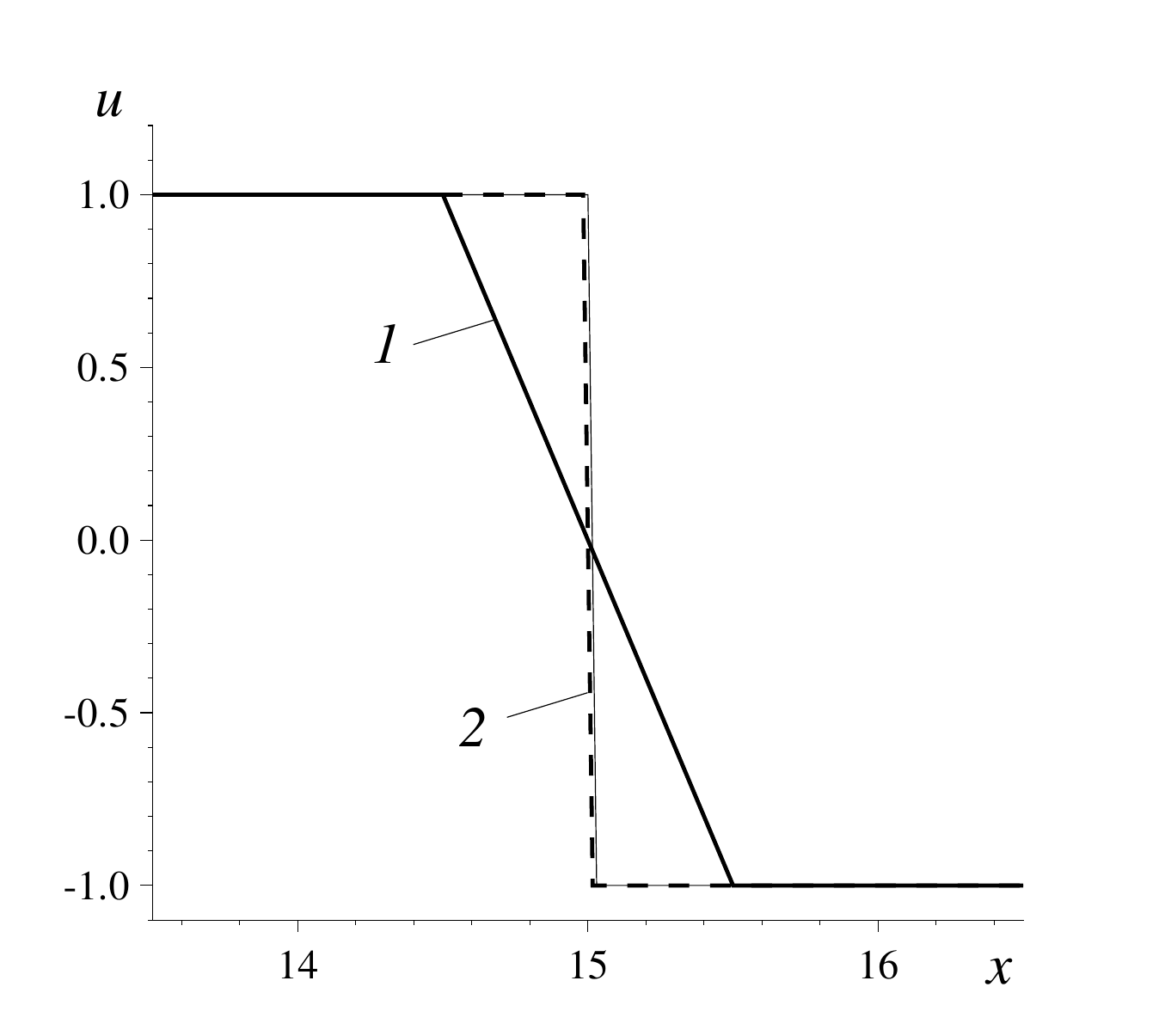}}
  \subfigure[]{\includegraphics[width=0.49\textwidth]{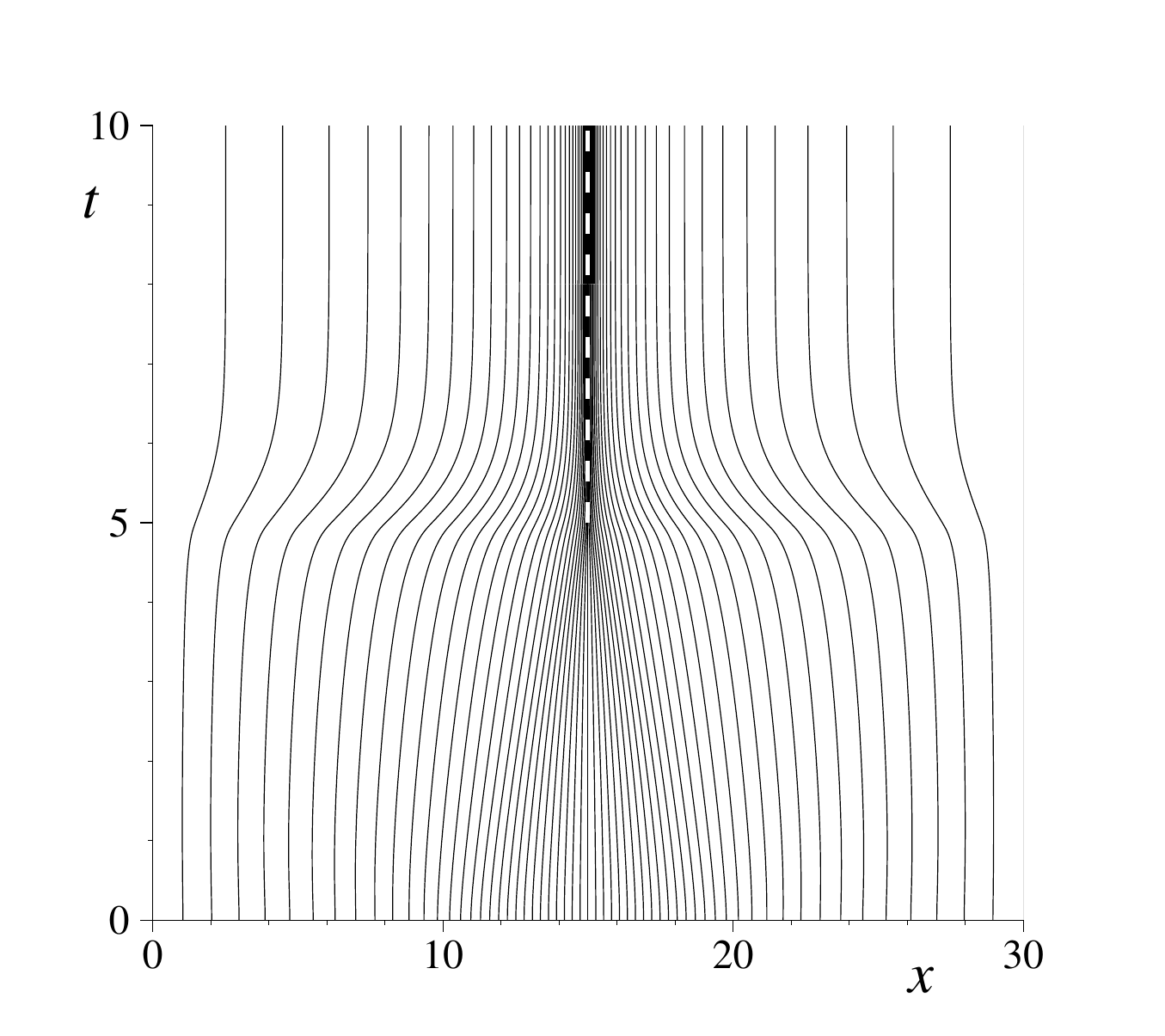}}
  \caption{\small\em Stationary shock wave formation in the \textsc{Burgers}--\textsc{Hopf} equation: (a) comparison of different profiles (thin solid line is the exact solution, (1) is the numerical solution on the uniform mesh, (2) is the solution on the moving mesh); (b): trajectories of grid nodes in space-time. Note that the panel (a) shows only a zoom of the computational domain.}
  \label{fig:3}
\end{figure}

In order to illustrate the performance of the generalized \textsc{Harten} scheme on moving meshes we take the classical example of the inviscid \textsc{Burgers}--\textsc{Hopf} equation ($f(u) = \frac{u^2}{2}$). As the first illustration let us take the following initial condition on the computational domain $x \in [0, \ell]$:
\begin{equation*}
  u_0(x)\ =\ \left\{
  \begin{array}{cl}
    u_l, & x\ \leq\ x_l, \\
    \displaystyle{\frac{x - x_r}{x_l - x_r}\; u_l\ +\ \frac{x - x_l}{x_r - x_l}\; u_r}, &  x_l\ <\ x\ <\ x_r, \\
    u_r & x\ \geq\ x_r.
  \end{array}
  \right.
\end{equation*}
All numerical parameters are given in Table~\ref{tab:params3}. The exact \emph{weak} solution to this problem is given by the following formula:
\begin{equation*}
  u(x,t)\ =\ \left\{
  \begin{array}{cl}
    u_l, & x\ \leq\ x_l\ +\ u_l\, t, \\
    \displaystyle{\frac{x\ -\ x_*}{t\ -\ t_*}}, & x_l\ +\ u_l\,t\ <\ x\ <\ x_r\ +\ u_r\, t, \\
    u_r, & x\ \geq\ x_r\ +\ u_r\,t,
  \end{array}
  \right.
  \qquad t\ <\ t_*.
\end{equation*}
The gradient catastrophe takes place at $t = t^*$ and $x = x_*$, which can be computed exactly as well
\begin{equation*}
  t^*\ =\ -\;\frac{x_r\ -\ x_l}{u_r\ -\ u_l}, \qquad
  x_*\ =\ x_l\ +\ u_l\,t^*\ \equiv\ x_r\ +\ u_r\, t^*.
\end{equation*}
After $t = t^*$ we obtain a travelling shock-wave profile:
\begin{equation*}
  u(x,t)\ =\ \left\{
  \begin{array}{cl}
    u_l, & x\ \leq\ x_*\ +\ \{u\}(t\ -\ t^*), \\
    u_r, & x\ >\ x_*\ +\ \{u\}(t\ -\ t^*),
  \end{array}
  \right.
  \qquad \{u\}\ \coloneqq\ \frac{u_l + u_r}{2}.
\end{equation*}
For the parameters chosen in our numerical simulation (see Table~\ref{tab:params3}) the shock wave remains stationary at $t\ =\ t^*\ =\ 5$ and $x\ =\ x_*\ =\ 15$ since $u_l\ =\ -u_r\,$. The numerical solution and the grid motion are shown in Figure~\ref{fig:3}. It can be seen that the solution monotonicity is preserved on both fixed and moving meshes. However, on the moving grid we obtain a solution which cannot be distinguished from the exact one up to the graphical resolution. On the fixed grid the shock wave has a finite width equal to two control volumes.

\begin{figure}
  \centering
  \subfigure[]{\includegraphics[width=0.49\textwidth]{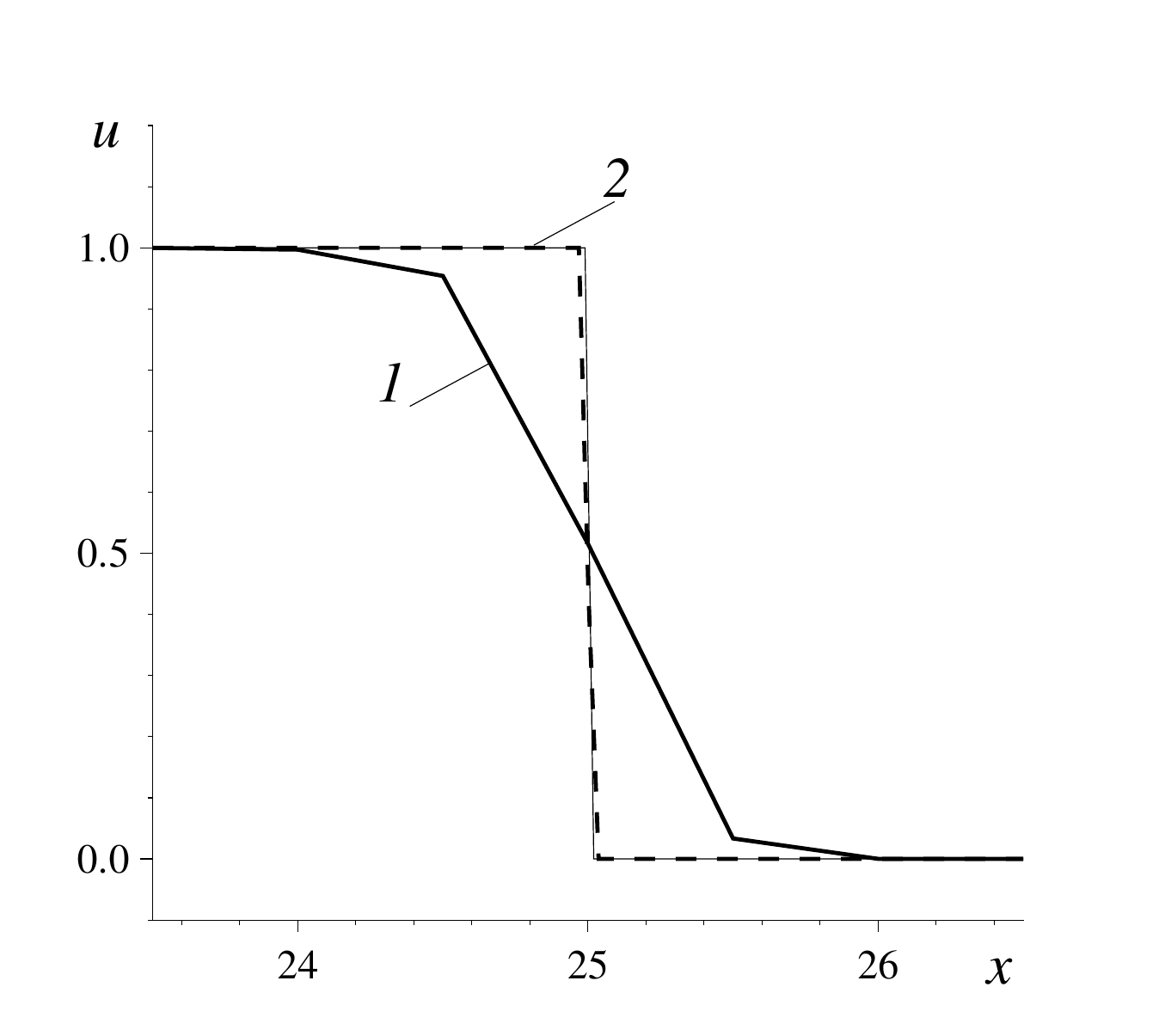}}
  \subfigure[]{\includegraphics[width=0.49\textwidth]{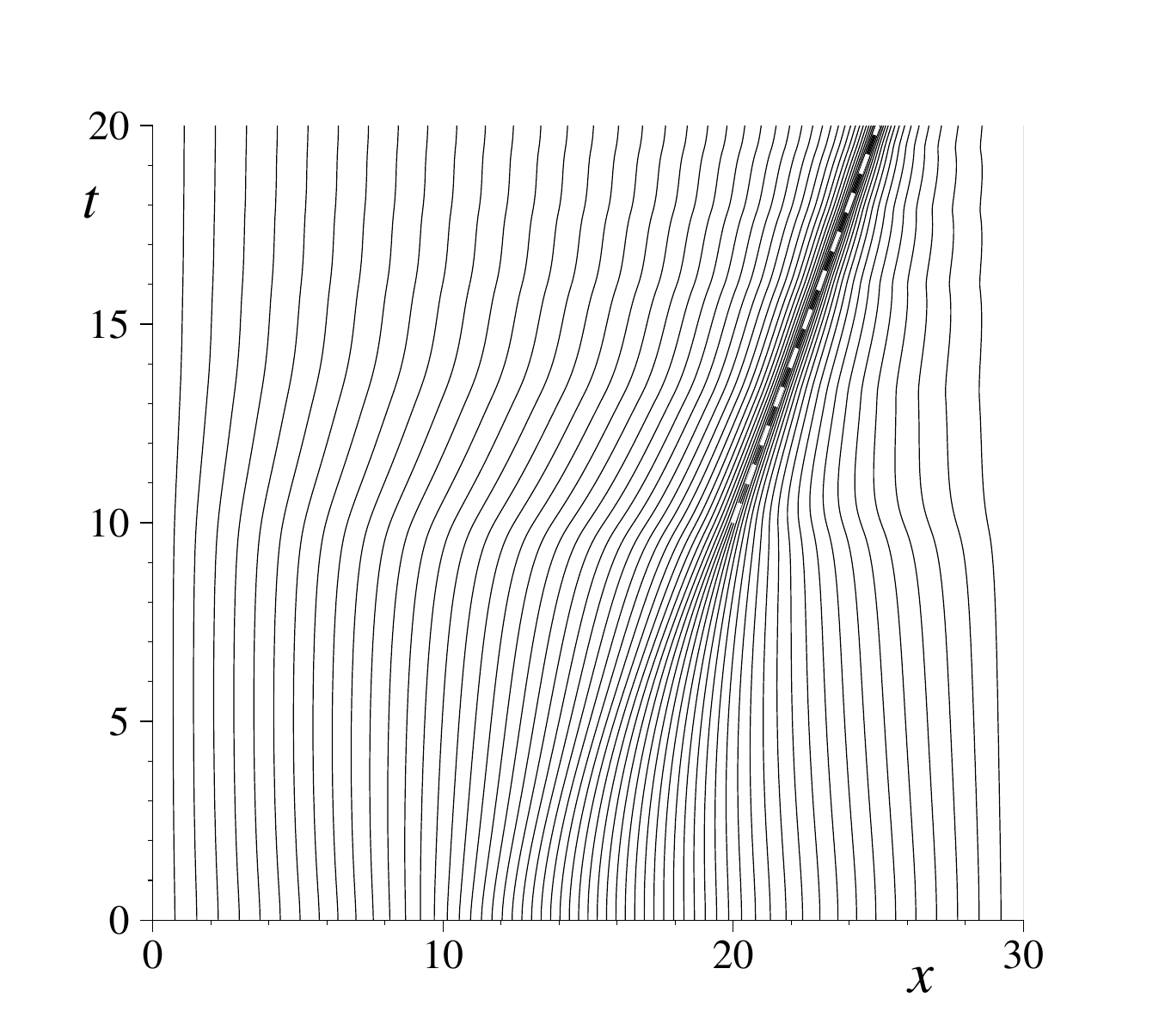}}
  \caption{\small\em Stationary shock wave formation in the \textsc{Burgers}--\textsc{Hopf} equation: (a) comparison of different profiles (thin solid line is the exact solution, (1) is the numerical solution on the uniform mesh, (2) is the solution on the moving mesh); (b): trajectories of grid nodes in space-time. Note that the panel (a) shows only a zoom of the computational domain.}
  \label{fig:4}
\end{figure}

As another numerical illustration, we consider the same problem, but we change the value of $u_r \coloneqq 0$ instead of $-1$. In this case the gradient catastrophe takes place at $t\ =\ t^*\ =\ 10$ and $x\ =\ x_*\ =\ 20\,$. The shock wave profile is not steady anymore since $\{u\} = \dfrac{u_l\ +\ u_r}{2}\ =\ \dfrac{1}{2}\ \neq\ 0\,$. So, it will move in the rightward direction with a constant speed equal to $\dfrac{1}{2}\,$. The results of numerical simulations are shown on Figure~\ref{fig:4}. It can be seen again that the adaptive solution is excellent, while the fixed grid produces a shock wave smeared out on four cells only.

\begin{table}
  \centering
  \begin{tabular}{lc}
  \hline\hline
  \textit{Parameter} & \textit{Value} \\
  \hline
  Computational domain length, $\ell$ & $30.0$ \\
  Number of nodes, $N$ & $60$ \\
  Final simulation time, $T_f$ & $10.0$ \\
  CFL number, $C$ & $0.2$ \\
  Left transition coordinate, $x_l$ & $10.0$ \\
  Right transition coordinate, $x_r$ & $20.0$ \\
  Left constant state, $u_l$ & $1.0$ \\
  Right constant state, $u_r$ & $-1.0$ \\
  Monitor function parameter, $\alpha$ & $15.0$ \\
  Grid diffusion parameter, $\beta$ & $80.0$ \\
  Smoothing parameter, $\sigma$ & $60.0$ \\
  \hline\hline
  \end{tabular}
  \bigskip
  \caption{\small\em Numerical parameters used to simulate a stationary shock-wave formation under the \textsc{Burgers}--\textsc{Hopf} equation dynamics.}
  \label{tab:params3}
\end{table}

\section{Nonlinear system of equations}
\label{sec:nswe}

The predictor-corrector scheme for nonlinear systems on moving meshes will be presented on a practically important case of the Nonlinear Shallow Water Equations (NSWE) \cite{Stoker1958b}. This choice reflects also the scientific interests of the authors of the present article. Moving grid techniques have been applied to other systems of conservation laws such as the \textsc{Euler} equations in \cite{Stockie2001}, pressureless gas dynamics \cite{Krejic2008} and magnetohydrodynamics \cite{VanDam2006}.

\begin{figure}
  \centering
  \includegraphics[width=0.95\textwidth]{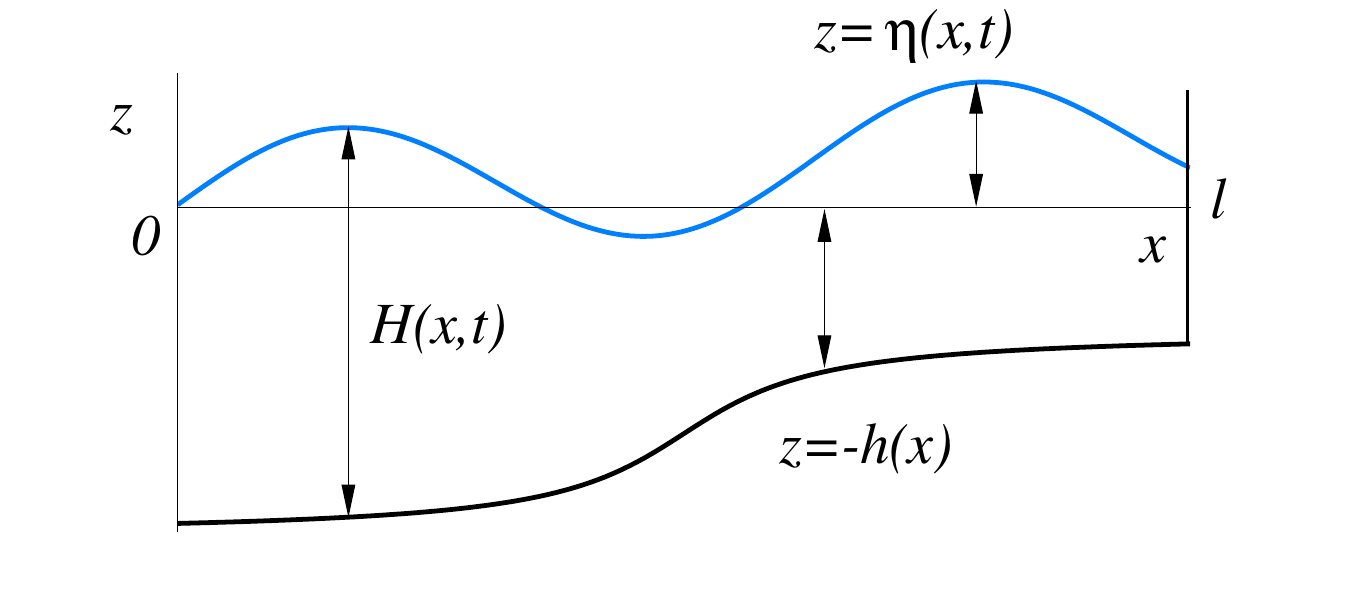}
  \caption{\small\em Sketch of the fluid domain in shallow water flows.}
  \label{fig:sketch}
\end{figure}

So, consider one-dimensional nonlinear shallow water equations written in a vectorial form:
\begin{equation}\label{eq:nswe}
  \u_t\ +\ \bigl[\,\f(\u)\,\bigr]_x\ =\ \G(\u, x)\,,
\end{equation}
where $\u$ is the vector of conservative variables, $\f(\u)$ is the flux function and $\G\,(\u,\,x)$ contains the source terms:
\begin{equation*}
  \u\ =\ \begin{pmatrix}
   H \\
   H u
  \end{pmatrix}, \qquad
  \f(\u)\ =\ \begin{pmatrix}
    H u \\
    H u^2\ +\ \frac{g}{2}H^2
  \end{pmatrix}, \qquad
  \G(\u, x)\ =\ \begin{pmatrix}
    0 \\
    g H h_x
  \end{pmatrix}\,.
\end{equation*}
Here $u(x,t)$ is the depth-averaged velocity and $H(x,t) = \eta(x,t) + h(x)$ is the total water depth defined as the sum of the bathymetry function $h(x)$ and the free surface elevation $\eta(x,t)$ over the mean water level. The sketch of the physical domain with all these definitions is shown in Figure~\ref{fig:sketch}. The system \eqref{eq:nswe} is considered on a finite segment $\I = [0, \ell]$ with appropriate boundary conditions. As before, we consider a smooth bijective time-dependent transformation $x = x(q,t)$ from the reference domain $Q = [0,1]$ into $\I$. The NSWE on the reference domain are written in the form \cite{Barakhnin2000}:
\begin{equation*}
  \v_t\ +\ \frac{1}{J}\;[\f_q\ -\ x_t\,\v_q]\ =\ \frac{1}{J}\;\G(\v,q),
\end{equation*}
where $\v(q,t) = \u\circ x = \u\bigl(x(q,t), t\bigr)$ and the source term contains a derivative with respect to $q\,$:
\begin{equation*}
  \G(\v, q)\ =\ \begin{pmatrix}
    0 \\
    g H h_q
  \end{pmatrix}, \qquad
  h(q,t)\ =\ h\bigl(x(q,t)\bigr).
\end{equation*}
It is noted that the bathymetry becomes time-dependent under this change of variables, even if initially the bottom was fixed. There exists also a conservative form NSWE on the reference domain:
\begin{equation*}
  (J\v)_t\ +\ [\f\ -\ x_t\,\v]_q\ =\ \G(\v,q), \qquad J\ \coloneqq\ x_q\ >\ 0.
\end{equation*}
Finally, the flux vector satisfies the following equation:
\begin{equation*}
  \f_t\ +\ \frac{1}{J}\;\A[\f_q\ -\ x_t\v_q]\ =\ \frac{1}{J}\;\A\G(\v,q),
\end{equation*}
where $\A$ is the advective flux's \textsc{Jacobian} matrix, which can be easily computed:
\begin{equation*}
  \A(\v)\ =\ \OD{\f(\v)}{\v}\ =\ \begin{pmatrix}
    0 & 1 \\
    -v^2 + gH & 2 v
  \end{pmatrix}.
\end{equation*}
Since the system of NSWE \eqref{eq:nswe} is hyperbolic, the \textsc{Jacobian} $\A$ can be decomposed into the product of left $\Ll$ and right $\Rr$ eigenvectors:
\begin{equation*}
  \A\ =\ \Rr\cdot\Dd\cdot\Ll,
\end{equation*}
where the matrices $\Ll$, $\Rr$ and $\Dd$ are defined as
\begin{equation*}
  \Ll\ =\ \frac{1}{c^2}\;\begin{pmatrix}
    -\lambda_2 & 1 \\
    -\lambda_1 & 1
  \end{pmatrix},
  \qquad
  \Rr\ =\ \frac{c}{2}\;\begin{pmatrix}
    -1 & 1 \\
    -\lambda_1 & \lambda_2
  \end{pmatrix}, \qquad
  \Dd\ =\ \begin{pmatrix}
    \lambda_1 & 0 \\
    0 & \lambda_2
  \end{pmatrix}.
\end{equation*}
Here $\lambda_{1,2} = v \mp c$ ($c = \sqrt{gH}$ being the speed of gravity waves) are eigenvalues of the \textsc{Jacobian} $\A\,(\v)\,$. They are real and distinct provided that $H(q,t) > 0$. Mathematically it implies the strict hyperbolicity property for the system \eqref{eq:nswe}.

\subsection{Finite volume discretization}

The reference domain $Q$ is discretized into $N$ equal control volumes $\bigl[\,q_{\,j},\, q_{\,j\,+\,1}\,\bigr]$ of the width $\Delta q\ \equiv\ \frac{1}{N}\,$. The nodes $x_{\,j}^{\,n}$ are obtained as images of uniform nodes under the mapping $x\,(q,\,t)\,$:
\begin{equation*}
  x_{\,j}^{\,n}\ =\ x\,(q_{\,j},\,t^{\,n})\,.
\end{equation*}
In the predictor-corrector scheme on the first step we evaluate the flux vector in cell centers:
\begin{equation*}
  \hat{\f}_{\,j\,+\,1/2}\ =\ \frac{\f^n_j\ +\ \f^n_{j+1}}{2}\ -\ \frac{\tau}{2}\;\Bigl[\frac{1}{J}\,\Rr\cdot\D\cdot{\bar \Dd}\cdot({\bar \Dd}\cdot\Pp\ -\ \Ll\cdot\G)\Bigr]^n_{j+1/2}, \quad \Pp\ \coloneqq\ \Ll\cdot\v_{\,q}\,.
\end{equation*}
Vector $\Pp$ is an analogue\footnote{The analogy is understood in the following sense: if we linearize our problem, the components of vector $\Pp$ will provide \emph{exactly} two \textsc{Riemann}'s invariants.} of the \textsc{Riemann}'s invariants in the nonlinear case. The corrector step is based on the conservative form of the equation:
\begin{equation}\label{eq:scheme}
  \frac{(J\,\v)_{\,j}^{\,n\,+\,1}\ -\ (J\,\v)_{\,j}^{\,n}}{\tau}\ +\ \frac{\bigl({\hat{\f}}\ -\ (x_{\,t}\,\v)^{\,n}\bigr)_{\,j\,+\,1/2}\ -\ \bigl({\hat{\f}}\ -\ (x_{\,t}\,\v)^{\,n}\bigr)_{\,j\,-\,1/2}}{\Delta q}\ =\ \G^{\,*}_{\,j}\,,
\end{equation}
where $\tau$ is the time step and the following discretizations are used 
\begin{equation*}
  \G_{\,j\,+\,1/2}^{\,n}\ =\ \begin{pmatrix}
  0 \\
  g\,(H\,h_{\,q})_{\,j\,+\,1/2}^{\,n}
  \end{pmatrix}, \quad
  \G_j^*\ =\ \begin{pmatrix}
  0 \\
  g(Hh_q)_j^{n+1/2}
  \end{pmatrix}\,, \quad
  \D_{j+1/2}^n\ =\ \begin{pmatrix}
    1 + \theta_{j+1/2}^{1,n} & 0 \\
    0 & 1 + \theta_{j+1/2}^{2,n}
  \end{pmatrix}\,.
\end{equation*}
The total water depth $H$ and local depth gradient $h_{\,q}$ have to be computed in the following way:
\begin{equation*}
  H^{n+1/2}_j = \frac{H^{n+1}_{j+1} + H^{n+1}_{j-1} + H^{n}_{j+1} + H^{n}_{j-1}}{4}, \ 
  h^{n+1/2}_{q,j} = \frac{h^{n+1}_{j+1} - h^{n+1}_{j-1} + h^{n}_{j+1} - h^{n}_{j-1}}{4\Delta q}, \ 
  h_{q,\, j+1/2}^n = \frac{h_{j+1}^n - h_j^n}{\Delta q}.
\end{equation*}
Despite the fact that $H_{j \pm 1}^{n+1}$ is present in the formulas above, the scheme remains explicit, since we compute first $\{H_j^{n+1}\}_{j=1}^{N-1}$ from the mass conservation and, then, this value is used in the correction of the momentum conservation equation.

The matrix ${\bar \Dd}$ plays the r\^ole of $\ab$ (the wave speed relative to the grid nodes) in the scalar case:
\begin{equation*}
  {\bar \Dd}\ \coloneqq\ \Dd_{j+1/2}^n\ -\ x_{t,\; j+1/2}^n\cdot\Id, \qquad
  \Dd_{j+1/2}^n\ =\ \begin{pmatrix}
    \lambda_{1,\; j+1/2}^n & 0 \\
    0 & \lambda_{2,\; j+1/2}^n
  \end{pmatrix}, \qquad
  \Id\ =\ \begin{pmatrix}
    1 & 0 \\
    0 & 1
  \end{pmatrix}.
\end{equation*}
Here $\lambda_{k,\; j+1/2}^n$, $k = 1,2$ are eigenvalues of the \textsc{Jacobian} matrix $\A_{j+1/2}^n$:
\begin{equation*}
  \A_{j+1/2}^n\ =\ \begin{pmatrix}
   0 & 1 \\
   -v_j^n v_{j+1}^n\ +\ g H_{j+1/2}^n & 2 v_{j+1/2}^n
  \end{pmatrix}, \qquad
  (\lambda_{1,2})_{j+1/2}^n\ =\ (v \mp c)_{j+1/2}^n\,,
\end{equation*}
the local gravity wave speed $c_{j+1/2}^n$ can be computed analytically:
\begin{equation*}
  c^n_{j+1/2}\ =\ \sqrt{(v^n_{j+1/2})^2\ -\ v^n_{j} v^n_{j+1}\ +\ gH^n_{j+1/2}}\ \geq\ \sqrt{gH^n_{j+1/2}}\ >\ 0\,.
\end{equation*}
This choice of the \textsc{Jacobian} $\A$ discretization is dictated by the requirement to preserve the compatibility condition $\f_{\,q}\ =\ \A\,(\v)\,\v_{\,q}$ at the discrete level as well. It is analogous to a similar algebraic condition on the \textsc{Jacobian} in the \textsc{Roe} scheme \cite{Roe1981}, but substantially different in details.

Finally, the scheme parameter $\theta_{j+1/2}^{k}$ is computed according to the same strategy as above, with the only difference being that it is computed separately for each component $k=1,\,2\,$:
\begin{equation}\label{eq:choice}
  \theta_{\,j\,+\,1/2}^{\,k}\ =\ \left\{
  \begin{array}{cl}
    0, & \abs{\gt_{j+1/2}^k}\ \leq\ \abs{\gt_{j+1/2-s_k}^k} \quad \textrm{and} \quad \gt_{j+1/2}^k \cdot \gt_{j+1/2-s_k}^k\ \geq\ 0, \\
    \theta_{0,j+1/2}^k (1\ -\ \xi^k_{j+1/2}), & \abs{\gt_{j+1/2}^k}\ >\ \abs{\gt_{j+1/2-s_k}^k} \quad \textrm{and} \quad \gt_{j+1/2}^k \cdot \gt_{j+1/2-s_k}^k\ \geq\ 0, \\
    \theta_{0,j+1/2}^k, & \gt_{j+1/2}^k \cdot \gt_{j+1/2-s_k}^k\ <\ 0\,,
\end{array}\right.
\end{equation}
where $\xi^k_{j+1/2} \coloneqq \frac{\gt_{j+1/2-s_k}^k}{\gt_{j+1/2}^k}$, $s_k \equiv s^n_{k,j+1/2} \coloneqq \sign\bigl({\bar \lambda}^n_{k,j+1/2}\bigr)$. Here ${\bar \lambda}^n_{k,j+1/2} = (\lambda_k - x_t)^n_{j+1/2}$ are simply the diagonal elements of the matrix ${\bar \Dd}$. The coefficient $\theta_{0,j+1/2}^k$ is related to the local CFL number as
\begin{equation*}
  \theta_{0,j+1/2}^k\ \coloneqq\ \frac{1}{C_{k,\; j+1/2}^{\;n}}\ -\ 1, \qquad
  C_{k,\; j+1/2}^{\;n}\ \coloneqq\ \frac{\tau}{\Delta q}\;\biggl(\frac{\abs{{\bar \lambda}_k}}{J}\;\biggr)^n_{j+1/2}\,.
\end{equation*}
Finally, we have to specify how to compute the indicators
\begin{equation*}
  \gt_{j+1/2}^k\ \coloneqq\ \bigl(\abs{{\bar \lambda}_k}\;(1\ -\ C_k)\;\widetilde{p}^k\bigr)_{j+1/2}^n,
\end{equation*}
where $\widetilde{p}^{\,k}\,$, $k\ =\ 1,\,2$ are components of the following vector:
\begin{equation}\label{eq:pp}
  \widetilde{\Pp}_{\,j\,+\,1/2}^{\,n}\ \coloneqq\ \frac{1}{(c_{\,j\,+\,1/2}^{\,n})^{\,2}}\;\begin{pmatrix}
    -\,c\,\eta_{\,q}\ +\ H\,v_{\,q} \\
    c\,\eta_{\,q}\ +\ H\,v_{\,q}
  \end{pmatrix}_{\,j\,+\,1/2}^{\,n}\,.
\end{equation}
The fully discrete scheme is stable under the following CFL-type restriction \cite{Courant1928} on the time step $\tau$:
\begin{equation*}
  \max_{k,\; j}\ \bigl\{\,C_{k,\; j\,+\,1/2}^{\,n}\,\bigr\}\ \leq\ 1\,.
\end{equation*}

\begin{remark}
We would like to explain better the motivation behind the choice of the vector $\widetilde{\Pp}$ (\cf Equation~\eqref{eq:pp}) in our scheme. This vector appears in the choice of the scheme parameter \eqref{eq:choice}, whose r\^ole consists in ensuring the monotonicity of the \textsc{Riemann}'s invariants. First of all, let us have a look at the components of vector $\Pp\,$:
\begin{equation*}
  \Pp_{\,j\,+\,1/2}^{\,n}\ =\ \bigl(\Ll\cdot\v_{\,q}\bigr)_{\,j\,+\,1/2}^{\,n}\ =\ \underbrace{\frac{1}{(c_{\,j\,+\,1/2}^{\,n})^{\,2}}\;\begin{pmatrix}
    -\,c\,\eta_{\,q}\ +\ H\,v_{\,q} \\
    c\,\eta_{\,q}\ +\ H\,v_{\,q}
  \end{pmatrix}_{\,j\,+\,1/2}^{\,n}}_{\equiv\ \widetilde{\Pp}_{\,j\,+\,1/2}^{\,n}}\ +\ 
  \Bigl(\,\frac{h_{\,q}}{c}\,\Bigr)_{\,j\,+\,1/2}^{\,n}\;\begin{pmatrix}
   -\,1 \\
   1
  \end{pmatrix}\,.
\end{equation*}
We can see that vector $\Pp$ is a sum of two vectors: the first addend depends on the solution gradient $(\eta_{\,q},\,v_{\,q})\,$, the second one depends on the bathymetry gradient $h_{\,q}\,$. We would like that the scheme parameters $\theta_{\,j\,+\,1/2}^{\,k}$ depend on the solution and not on the bathymetry. Thus, we made a choice to include only the first addend in the vector $\widetilde{\Pp}\,$.
\end{remark}

\subsubsection{Well-balanced property}

The following property can be shown for the scheme presented above:
\begin{theorem}\label{thm:well}
Assume that we have a fixed, but possibly non-uniform, grid $\{x_j^n\}$. If on the $n$\up{th} time layer we have the \emph{lake at rest state}, \ie
\begin{equation}\label{eq:lake}
  \eta_j^n\ \equiv\ 0\,, \qquad v_j^n\ \equiv\ 0\,,
\end{equation}
then, the predictor--corrector scheme described above will preserve this state at the $(n+1)$\up{th} time layer.
\end{theorem}

\begin{proof}
Since $x_{t,j} \equiv 0$, $\forall j$ and taking into account the definition \eqref{eq:lake} of the lake at rest state, we have:
\begin{equation*}
  \bigl({\bar \Dd} \cdot \Pp\ -\ \Ll \cdot \G\bigr)^n_{j+1/2}\ \equiv\ 0.
\end{equation*}
Consequently, the predictor step leads
\begin{equation*}
  {\hat\f}_{j+1/2}\ =\ \frac{\f^n_j\ +\ \f^n_{j+1}}{2}.
\end{equation*}
Then, from the obvious identities $h_j^{n+1} \equiv h_j^n$, $J_j^{n+1} \equiv J_j^n$ and the way how we compute $\G^*_j$, it follows that $\eta_{j}^{n+1} \equiv 0$, $v_{j}^{n+1} \equiv 0$. Thus, the lake at rest state is preserved at the next time layer as well.
\end{proof}
The well-balanced property is absolutely crucial to compute qualitatively correct numerical solutions to conservation (balance) laws \cite{Gosse2013}. The Theorem~\ref{thm:well} shows that the discretization proposed above does possess this property regardless some complications coming from the grid motion.

\subsubsection{Conservation property on moving grids}

For conservation laws with source terms (the so-called \emph{balance laws}), the conservation property of a numerical scheme is understood in the sense that the \emph{fully} discrete solution verifies some thoroughly chosen finite difference analogues of continuous balance laws. For the NSWE we have the conservation of mass and the balance of horizontal (depth-averaged) momentum. Please, note that on a flat bottom the balance of momentum becomes a conservation law as well. The conservation property is crucial as well to compute qualitatively correct numerical solutions \cite{Gosse2013}.

On moving grids we have an additional difficulty: the computation of the speeds of mesh nodes has to be consistent with the calculation of the cell volumes. If it is done in the right way, the discrete geometric conservation law will be verified as well \cite{Thomas1979}. In 1D it can be written as
\begin{equation}\label{eq:gcl}
  J_j^{n+1}\ =\ J_j^n\ +\ \tau\, x_{t q,j}^n.
\end{equation}

We shall show below that other schemes can by put into the conservative form similar to \eqref{eq:scheme}. The peculiarity here consists in the fact that the \textsc{Jacobian} matrix $\A$ is non-constant, in contrast to the linear advection equation, for which any reasonable scheme is conservative (see also Remark~\ref{rem:cons} below about a nonlinear scalar equation). Moreover, the moving mesh is another complication, which makes this task rather non-trivial. For example, the scheme parameter $\theta_{j\pm 1/2}^k$ can be chosen as
\begin{equation*}
  \theta_{j \pm 1/2}^k\ =\ \theta_{0, j \pm 1/2}^k, \qquad k\ =\ 1,2.
\end{equation*}
This choice corresponds to the first order upwind scheme (as it was discussed above in Section~\ref{sec:predcorr}). So, the scheme \eqref{eq:scheme} can be rewritten as
\begin{multline*}
  \frac{(J\v)_{j}^{n+1}\ -\ (J\v)_{j}^n}{\tau}\ +\\ 
  \frac{1}{h}\;\biggl[\Bigl(\frac{\f^n_{j+1} + \f^n_{j}}{2}\ -\ \frac{\bigl({\bar{\A}}^+ - {\bar{\A}}^-\bigr)^n_{j+1/2}}{2}(\v^n_{j+1} - \v^n_{j})\ -\ (x_t\,\v)^n_{j+1/2}\ +\ \frac{h}{2}(\Rr\S\Ll\G)^n_{j+1/2}\Bigr)\ -\\
  \Bigl(\frac{\f^n_{j} + \f^n_{j-1}}{2}\ -\ \frac{\bigl({\bar{\A}}^+ - {\bar{\A}}^-\bigr)^n_{j-1/2}}{2}(\v^n_{j} - \v^n_{j-1})\ -\ (x_t\,\v)^n_{j-1/2}\ +\ \frac{h}{2}(\Rr\S\Ll\G)^n_{j-1/2}\Bigr)\biggr] = \G_j^*,
\end{multline*}
where ${\bar{\A}}^\pm = \Rr\bar{\Lambda}^\pm\Ll$, $\bar{\lambda}_k^{\pm} = \half(\bar{\lambda}_k \pm \abs{\bar{\lambda}_k})$ and
\begin{equation*}
  \S\ =\ \begin{pmatrix}
    \sign(\lambda_1) & 0 \\
    0 & \sign(\lambda_2)
  \end{pmatrix}, \qquad
  \bar{\Lambda}^\pm\ =\ \begin{pmatrix}
    \bar{\lambda}_1^\pm & 0 \\
    0 & \bar{\lambda}_2^\pm
  \end{pmatrix}.
\end{equation*}
One can see also that the following identities hold:
\begin{equation*}
  \bar{\Lambda}\ \equiv\ \bar{\Lambda}^+\ +\ \bar{\Lambda}^-, \qquad
  \bar{\A}^+\ +\ \bar{\A}^-\ \equiv\ \bar{\A}\ \equiv\ \A\ -\ x_t\,\Id.
\end{equation*}
By using these identities and the discrete geometric conservation law \eqref{eq:gcl}, one can derive also the following equivalent form of the scheme \eqref{eq:scheme}:
\begin{multline*}
  \frac{\v_{j}^{n+1}\ -\ \v_{j}^n}{\tau}\ +\ \frac{1}{J^{n+1}_j}\;\Bigl[\bigl(\A^{+}\v_q\bigr)^n_{j-1/2}\ +\ \bigl(\A^{-}\v_q\bigr)^n_{j+1/2}\Bigr]\ = \\
  \frac{1}{J^{n+1}_j}\;\Bigl[\G^*_j\ -\ \half\bigl((\Rr\S\Ll\G)^n_{j+1/2}\ -\ (\Rr\S\Ll\G)^n_{j-1/2}\bigr)\Bigr].
\end{multline*}

The choice \eqref{eq:choice} of the scheme parameter $\theta_{j\pm 1/2}^k$ leads to a conservative scheme as well, but the underlying computations (to show it) are less elegant.

\begin{remark}\label{rem:cons}
In order to illustrate even better the interplay between conservative and non-conservative forms, let us take again the inviscid \textsc{Burgers}--\textsc{Hopf} equation \cite{LeVeque1992}:
\begin{equation*}
  u_t\ +\ \Bigl(\frac{u^2}{2}\Bigr)_x\ =\ 0.
\end{equation*}
For the sake of simplicity let us discretize this equation on a uniform grid using the following non-conservative scheme:
\begin{equation*}
  \frac{u_j^{n+1}\ -\ u_j^n}{\tau}\ +\ \frac{u_j^n\ +\ u_{j-1}^n}{2}\cdot 
  \frac{u_j^n\ -\ u_{j-1}^n}{\Delta x}\ =\ 0.
\end{equation*}
It is not difficult to see that this non-conservative scheme is equivalent to the following conservative upwind scheme:
\begin{equation*}
  \frac{u_j^{n+1}\ -\ u_j^n}{\tau}\ +\ \frac{1}{\Delta x}\;\biggl(\,\frac{(u_j^n)^2}{2}\ -\ \frac{(u_{j-1}^n)^2}{2}\,\biggr)\ =\ 0.
\end{equation*}
However, if one takes a different non-conservative scheme\footnote{The difference is in the estimation of propagation speed, \cf $u_j^n$ \vs $\half(u_j^n + u_{j-1}^n)$.}
\begin{equation}\label{eq:bad}
  \frac{u_j^{n+1}\ -\ u_j^n}{\tau}\ +\ u_j^n\cdot \frac{u_j^n\ -\ u_{j-1}^n}{\Delta x}\ =\ 0,
\end{equation}
it is straightforward to see that one cannot transform it in a conservative form. Moreover, one can even show that scheme \eqref{eq:bad} is not converging. So, the conservation and convergence properties are intrinsically related.
\end{remark}

\subsection{Numerical results}

In this section we present two test cases of complex wave propagation in shallow water equations. The first test case is a validation against an implicit analytic solution, derived below specifically for this purpose, while the second test case is related to real world applications.

\subsubsection{Exact solution derivation}
\label{sec:exa}

In order to validate the numerical scheme, we derive an exact solution to nonlinear shallow water equations which will serve as a reference solution. We already saw this particular solution in the literature \cite{Voltsinger1989}, however, the derivation procedure is published for the first time to our knowledge.

First of all, we assume that the bottom is flat, \ie $h(x) = h_0 = \const > 0$ and we consider an initial value problem on the real line:
\begin{equation*}
  \eta\,(x,\,0)\ =\ \eta_{\,0}\,(x)\,, \qquad u\,(x,\,0)\ =\ u_{\,0}\,(x)\,.
\end{equation*}
The nonlinear shallow water equations are rewritten using the \textsc{Riemann} invariants \cite{Toro2009}:
\begin{equation*}
  r_{\,t}\ +\ \frac{3\,r\ +\ s}{4}\; r_{\,x}\ =\ 0\,, \qquad
  s_{\,t}\ +\ \frac{r\ +\ 3\,s}{4}\; s_{\,x}\ =\ 0\,,
\end{equation*}
where 
\begin{equation*}
  r\ \coloneqq\ u\ -\ 2\,c\,, \qquad
  s\ \coloneqq\ u\ +\ 2\,c\,, \qquad
  c\ \coloneqq\ \sqrt{g\,H}\,.
\end{equation*}
We seek for a particular solution in the form of an $r$-wave propagating in the leftward direction. In other words we assume $s = s_0 = \const$. Assuming that at the infinity the water is at rest, we deduce that $s_0 = 2c_0 = 2\sqrt{g h_0}$. Then, from formula $u + 2c = s_0$ we deduce that the speed $u(x,t)$ is related to the total depth:
\begin{equation}\label{eq:speed}
  u(x,t)\ =\ s_0\ -\ 2 c\ =\ 2c_0\ -\ 2\sqrt{gH(x,t)}.
\end{equation}
Consequently, taking the limit $t\to 0$ we have a similar relation between the initial conditions:
\begin{equation*}
  u_0(x)\ =\ 2c_0\ -\ 2\sqrt{g(\eta_0(x) + h_0)}.
\end{equation*}
If we introduce a new variable
\begin{equation*}
  p(x,t)\ \coloneqq\ \frac{3 r(x,t)\ +\ s_0}{4}
\end{equation*}
the governing equation for the $r$-invariant becomes simply the \textsc{Burgers}--\textsc{Hopf} equation:
\begin{equation*}
  p_t\ +\ p\,p_x\ =\ 0.
\end{equation*}
The solution can be obtained using the method of characteristics by the following implicit formula \cite{Rozhdestvenskiy1978, Whitham1999}:
\begin{equation*}
  p(x,t)\ =\ p_0\bigl(x\ -\ p(x,t)t\bigr),
\end{equation*}
where $p_0(x)$ is the initial condition
\begin{equation*}
  p_0(x)\ =\ \frac{3 r(x,0)\ +\ s_0}{4}\ =\ 2c_0\ -\ 3\sqrt{g(\eta_0(x) + h_0)}.
\end{equation*}
Consequently, the solution $p(x,t)$ can be obtained at any space location $x$ and time instance $t$ by solving the following implicit equation
\begin{equation}\label{eq:eq}
  p(x,t)\ =\ 2c_0\ -\ 3\sqrt{g(\eta_0\bigl(x - p(x,t)\,t\bigr) + h_0)}.
\end{equation}
Once, $p(x,t)$ is found, the free surface elevation $\eta(x,t)$ can be reconstructed by using this formula
\begin{equation}\label{eq:etaref}
  \eta(x,t)\ =\ \biggl(\frac{2c_0\ -\ p(x,t)}{3\sqrt{g}}\biggr)^2\ -\ h_0,
\end{equation}
and the flow speed $u(x,t)$ from \eqref{eq:speed}.

\begin{remark}
Please, note that this method works only until the gradient catastrophe occurs at some time instance $t = t^*$, which depends on the initial condition.
\end{remark}

\begin{remark}
In order to guarantee the solvability of equation \eqref{eq:eq}, it is enough to check that the following function
\begin{equation*}
  f(p; x, t)\ =\ p\ +\ 3\sqrt{g(\eta_0\bigl(x - p\,t\bigr) + h_0)}\ -\ 2c_0
\end{equation*}
is monotonic in $p$ for $0 < t < t^*$.
\end{remark}

\subsubsection{Validation}

\begin{figure}
  \centering
  \subfigure[]{\includegraphics[width=0.49\textwidth]{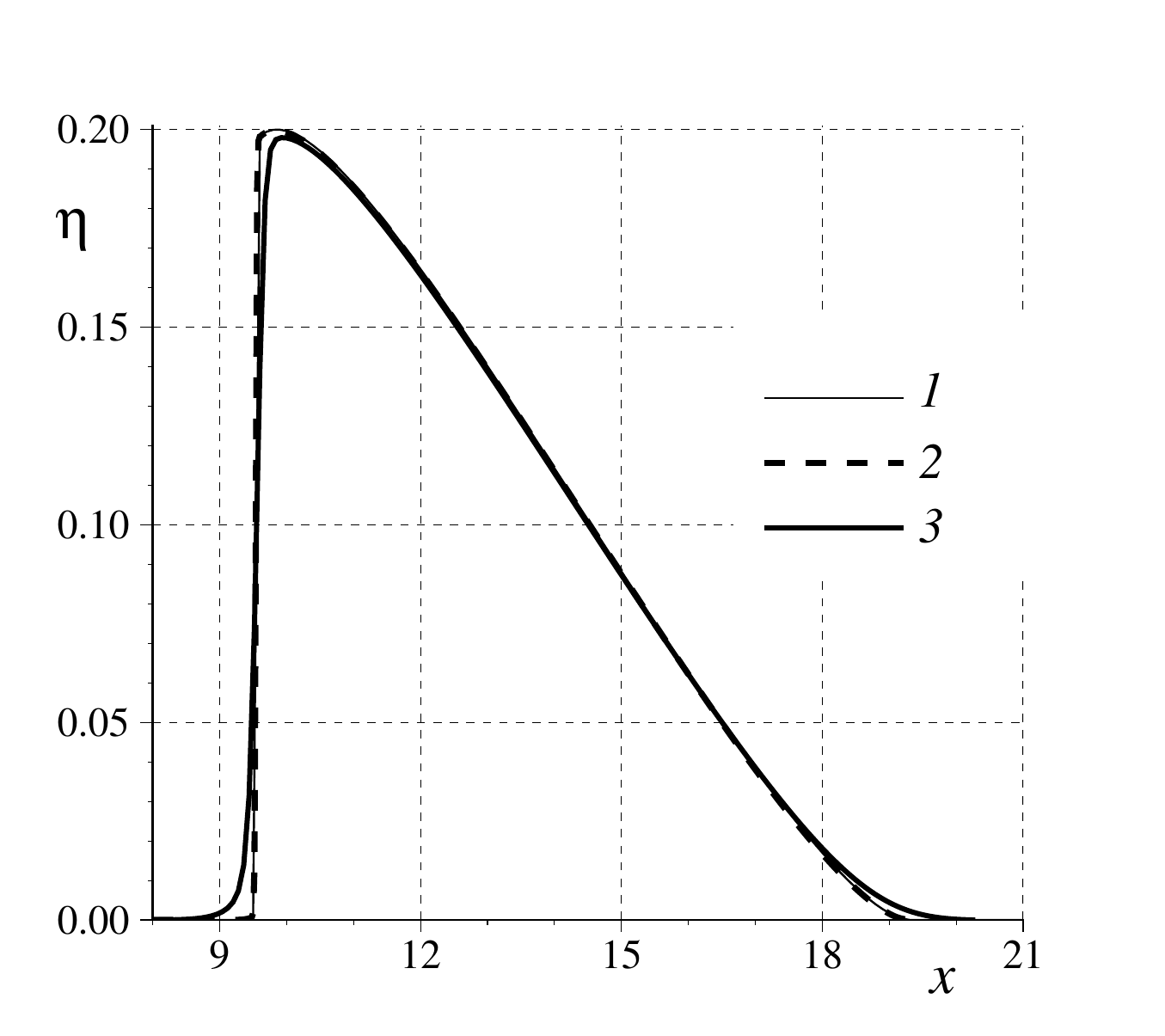}}
  \subfigure[]{\includegraphics[width=0.49\textwidth]{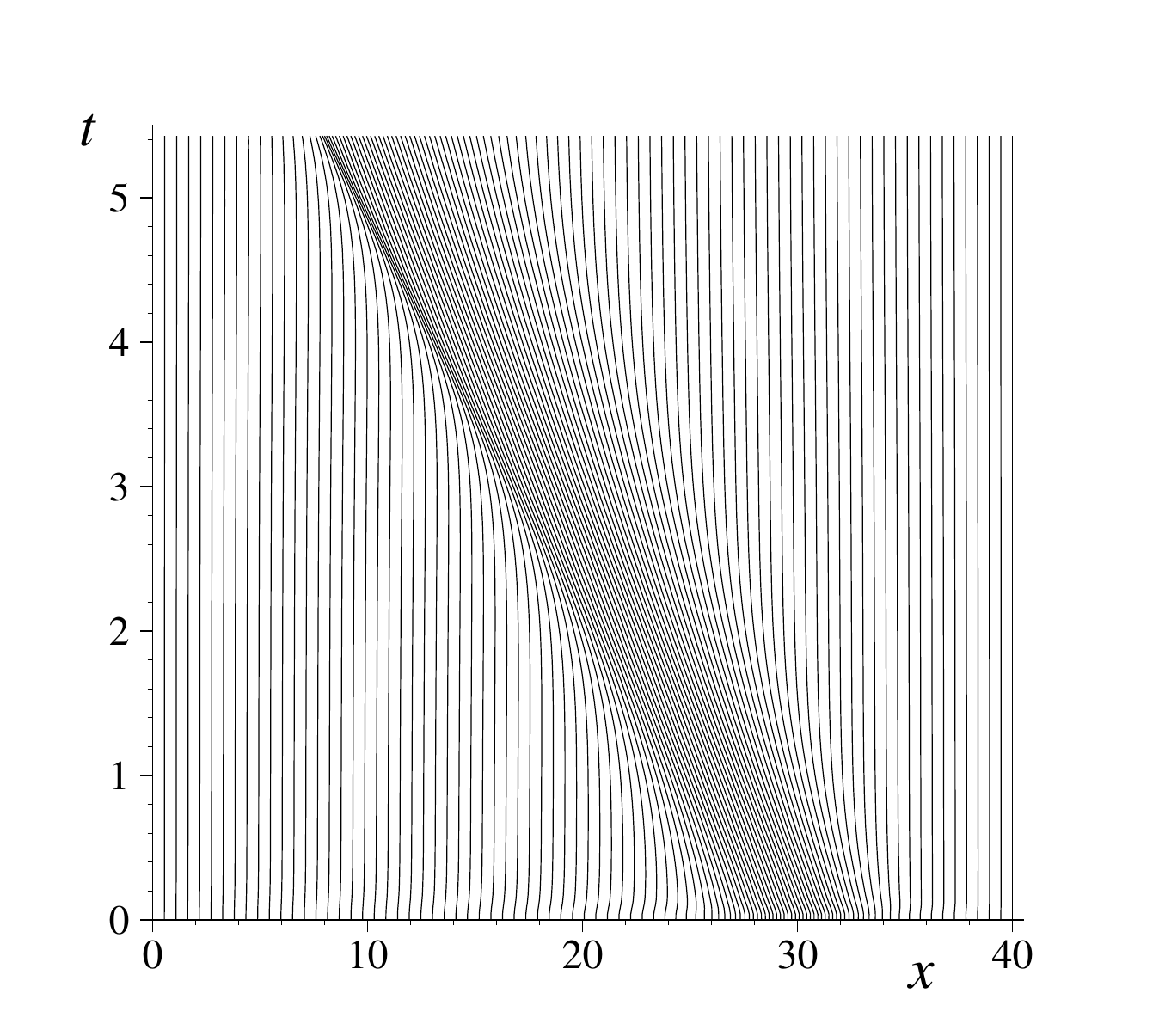}}
  \caption{\small\em Shock wave formation in nonlinear shallow water equations: (a) comparison of different profiles (thin solid line (1) is the exact solution, (2) is the numerical solution on the moving mesh, (3) is the solution on the fixed uniform mesh); (b): trajectories of grid nodes in space-time. Notice that the panel (a) shows only a zoom of the computational domain.}
  \label{fig:5}
\end{figure}

Using the method described above, we obtain the analytical solution for the following initial value problem:
\begin{equation*}
  \eta_0(x)\ =\ \left\{
    \begin{array}{lr}
    \frac{a}{2}\,\Bigl[1\ +\ \cos\bigl(\frac{2\pi(x\ -\ x_w)}{\lambda}\bigr)\Bigr],  & \abs{x\ -\ x_w}\ \leq\ {\lambda}/{2}, \\
    0, & \abs{x\ -\ x_w}\ >\ {\lambda}/{2},
  \end{array}
  \right.
\end{equation*}
where $x_w$ is the wave crest position and $a$ is the wave amplitude. For $0 < t < t^*$ the profile $\eta_0(x)$ will propagate to the left with the speed $c_0$ preserving its wavelength $\lambda$ and amplitude $a$. However, the gradient catastrophe will inevitably take place since the characteristics $\od{x}{t} = p_0(x)$ departing from the interval $(x_w - \lambda/2,\; x_w)$ form a converging beam. On the other hand, the characteristics from $(x_w,\; x_w + \lambda/2)$ diverge. Thus, the solution looks like a compression wave followed by a rarefaction wave. With the parameters specified in Table~\ref{tab:params4} the gradient catastrophe takes place about $t = t^* \approx 5.57$. Numerical results are shown in Figure~\ref{fig:5} right before the gradient catastrophe occurs. The adaptive grid was constructed using the following monitor function, which contains an additional term with respect to \eqref{eq:def}
\begin{equation}\label{eq:dimeq1}
  w(x,t)\ =\ 1\ +\ \alpha_0\,\abs{\eta(x,t)}\ +\ \alpha_1\,\abs{\eta_x(x,t)},
\end{equation}
which is discretized as
\begin{equation}\label{eq:dimeq2}
  w^n_{j+1/2}\ =\ 1\ +\ \alpha_0\, \abs{\eta_{j+1/2}^{n}}\ +\ \alpha_1\, \frac{\abs{\eta_{j+1}^{n}\ -\ \eta_{j}^{n}}}{x_{j+1}^{n}\ -\ x_{j}^{n}}.
\end{equation}
Thin solid line in Figure~\ref{fig:5} represents the exact solution computed by solving equation \eqref{eq:eq} as it was explained in the previous Section~\ref{sec:exa}. The agreement between these two profiles is exemplary. On the right panel (\textit{b}) one can see how the nodes follow the wave propagation and how their density increases around the steep gradient. On the other hand, the use of a fixed grid leads to a noticeable smearing of the numerical solution.

\begin{table}
  \centering
  \begin{tabular}{lc}
  \hline\hline
  \textit{Parameter} & \textit{Value} \\
  \hline
  Computational domain length, $\ell$ & $40.0$ \\
  Number of nodes, $N$ & $100$ \\
  Final simulation time, $T_f$ & $5.0$ \\
  CFL number, $C$ & $0.95$ \\
  Gravity acceleration, $g$ & $9.81$ \\
  Initial wave amplitude, $a$ & $0.2$ \\
  Wave crest initial position, $x_w$ & $30.0$ \\
  Initial wavelength, $\lambda$ & $10.0$ \\
  Undisturbed water depth, $h_0$ & $1.0$ \\
  Monitor function parameter, $\alpha_0$ & $10.0$ \\
  Monitor function parameter, $\alpha_1$ & $10.0$ \\
  Grid diffusion parameter, $\beta$ & $5.0$ \\
  Smoothing parameter, $\sigma$ & $5.0$ \\
  \hline\hline
  \end{tabular}
  \bigskip
  \caption{\small\em Numerical parameters used in validations of the nonlinear shallow water equations solver.}
  \label{tab:params4}
\end{table}

The exact solution presented in the previous Section~\ref{sec:exa} can be used in order to study the absolute accuracy of proposed numerical methods. To this end we introduce two discrete norms:
\begin{align*}
  \norm{f_{\,h}\ -\ \f_{\,\mathrm{ref}}}_{\,\infty}\ &:=\ \max_{1\ \leq\ j\ \leq\ N\,-\,1} \abs{f_{\,j}^{\,n}\ -\ f\,(x_{\,j},\,T_{\,f})}\,, \\
  \norm{f_{\,h}\ -\ \f_{\,\mathrm{ref}}}_{\,2}\ &:=\ \biggl(\,\sum_{j\,=\,1}^{N\,-\,1}\,\bigl(f_{\,j}^{\,n}\ -\ f\,(x_{\,j},\,T_{\,f})\bigr)^{\,2}\cdot\frac{x_{\,j+1}^{\,n}\ -\ x_{\,j-1}^{\,n}}{2}\,\biggr)^{\,2}\,,
\end{align*}
where $n$ is the time step number corresponding to the final simulation time $t\ =\ T_{\,f}\,$. Here $f$ denotes the free surface elevation $\eta$ or the depth-averaged horizontal velocity $u\,$. The reference solution for $u$ and $\eta$ is provided by formulas \eqref{eq:speed} and \eqref{eq:etaref} correspondingly. We would like to mention that the above definition of the norm $\norm{\cdot}_{\,2}$ is valid for the fixed and moving grids. The computational results are reported for the fixed and moving grids in Tables~\ref{tab:fix} and \ref{tab:mov} respectively. All the parameters in these simulations were fixed and taken from Table~\ref{tab:params4} (except the parameter $N\,$, which was varied). The second order convergence can be observed in Table~\ref{tab:fix}. In general, the results obtained on a moving grid are more accurate than those on a fixed one as it follows from the comparison of corresponding cells in Tables~\ref{tab:fix} and \ref{tab:mov}. In the last column of these tables we report also the CPU times measured on a computer equiped with an Intel Pentium CPU \@ 3.2 Ghz and 1Gb of RAM. The code was implemented in Fortran and compiled with Compaq Visual Fortran compiler (version 6.6c). In particular, for the finest grid ($N\ =\ 6\,400$) the overhead related to the grid motion represents less than 7\% of the overall CPU time and the error (measured in the most stringent norm $\norm{\cdot}_{\,\infty}$) is reduced by a factor of four. Thus, we may conclude that the application of moving grid methods is very interesting from the practical point of view.

We underline the fact that even better results could be obtained if one optimizes the monitor function parameters $\alpha_{\,0,\,1}$ along with mesh motion parameters $\beta$ and $\sigma\,$. In our numerical simulations we just took some \emph{reasonable} values. At the current stage of the development of these methods, it is difficult to say something on the optimal values of these parameters \emph{a priori}. Finally, we would like to mention that this test case is relatively difficult for both analytical and numerical methods since we deal with a nonlinear system of hyperbolic equations and at the final time $t\ =\ T_{\,f}$ the wave has almost developed a discontinuity (from a smooth initial condition). Thus, we operate at the edge of validity of the analytical solution \eqref{eq:speed}, \eqref{eq:etaref}.

\begin{table}
\centering
\caption{\small\em Numerical solution error measured at the final simulation time $t\ =\ T_{\,f}$ on a fixed grid. The reference solution is given be formulas \eqref{eq:speed} and \eqref{eq:etaref}. All numerical parameters for this simulation are reported in Table~\ref{tab:params4}.}
\label{tab:fix}
\bigskip
\begin{adjustwidth}{-1.0cm}{}
{\tabcolsep=5mm
\begin{tabular}{cccccc}
\hline\hline
 $N$ & \multicolumn{2}{c}{$\eta$} &  \multicolumn{2}{c}{$u$} \\ \cline{2-5}
 & $\norm{\eta_{\,h}-\eta_{\,\mathrm{ref}}}_{\,\infty}$ & $\norm{\eta_{\,h}-\eta_{\,\mathrm{ref}}}_{\,2}$ & $\norm{u_{\,h}-u_{\,\mathrm{ref}}}_{\,\infty}$ & $\norm{u_{\,h}-u_{\,\mathrm{ref}}}_{\,2}$ & CPU time, \textsf{s} \\ \hline
  100 & $0.388\cdot 10^{-1}$  & $0.289\cdot 10^{-1}$ & $0.363\cdot 10^{-1}$ & $0.270\cdot 10^{-1}$ & $<\ 10^{\,-\,2}$ \\ \hline
  200 & $0.245\cdot 10^{-1}$  & $0.127\cdot 10^{-1}$ & $0.229\cdot 10^{-1}$ & $0.118\cdot 10^{-1}$ & $0.03$ \\ \hline
  400 & $0.104\cdot 10^{-1}$  & $0.451\cdot 10^{-2}$ & $0.977\cdot 10^{-2}$ & $0.424\cdot 10^{-2}$ & $0.20$ \\ \hline
  800 & $0.418\cdot 10^{-2}$  & $0.144\cdot 10^{-2}$ & $0.394\cdot 10^{-2}$ & $0.136\cdot 10^{-2}$ & $1.36$ \\ \hline
  1600 & $0.128\cdot 10^{-2}$ & $0.393\cdot 10^{-3}$ & $0.122\cdot 10^{-2}$ & $0.373\cdot 10^{-3}$ & $9.05$ \\ \hline
  3200 & $0.353\cdot 10^{-3}$ & $0.103\cdot 10^{-3}$ & $0.337\cdot 10^{-3}$ & $0.981\cdot 10^{-4}$ & $71.08$ \\ \hline
  6400 & $0.885\cdot 10^{-4}$ & $0.257\cdot 10^{-4}$ & $0.843\cdot 10^{-4}$ & $0.245\cdot 10^{-4}$ & $573.11$ \\
  \hline\hline
\end{tabular}}
\end{adjustwidth}
\bigskip\bigskip
\end{table}

\begin{table}
\centering
\caption{\small\em Numerical solution error measured at the final simulation time $t\ =\ T_{\,f}$ on a moving grid. The reference solution is given be formulas \eqref{eq:speed} and \eqref{eq:etaref}. All numerical parameters for this simulation are reported in Table~\ref{tab:params4}.}
\label{tab:mov}
\bigskip
\begin{adjustwidth}{-1.0cm}{}
{\tabcolsep=5mm
\begin{tabular}{cccccc}
\hline\hline
 $N$ & \multicolumn{2}{c}{$\eta$} &  \multicolumn{2}{c}{$u$} \\ \cline{2-5}
 & $\norm{\eta_{\,h}-\eta_{\,\mathrm{ref}}}_{\,\infty}$ & $\norm{\eta_{\,h}-\eta_{\,\mathrm{ref}}}_{\,2}$ & $\norm{u_{\,h}-u_{\,\mathrm{ref}}}_{\,\infty}$ & $\norm{u_{\,h}-u_{\,\mathrm{ref}}}_{\,2}$ & CPU time, \textsf{s} \\ \hline
  100 & $0.125\cdot 10^{-1}$  & $0.907\cdot 10^{-2}$ & $0.123\cdot 10^{-1}$ & $0.900\cdot 10^{-2}$ & $0.03$ \\ \hline
  200 & $0.591\cdot 10^{-2}$ & $0.415\cdot 10^{-2}$ & $0.586\cdot 10^{-2}$ & $0.413\cdot 10^{-2}$ & $0.20$ \\ \hline
  400 & $0.266\cdot 10^{-2}$  & $0.170\cdot 10^{-2}$ & $0.264\cdot 10^{-2}$ & $0.170\cdot 10^{-2}$ & $1.13$ \\ \hline
  800 & $0.897\cdot 10^{-3}$  & $0.608\cdot 10^{-3}$ & $0.894\cdot 10^{-3}$ & $0.606\cdot 10^{-3}$ & $5.76$ \\ \hline
  1600 & $0.317\cdot 10^{-3}$ & $0.215\cdot 10^{-3}$ & $0.316\cdot 10^{-3}$ & $0.214\cdot 10^{-3}$ & $27.66$ \\ \hline
  3200 & $0.881\cdot 10^{-4}$ & $0.566\cdot 10^{-4}$ & $0.878\cdot 10^{-4}$ & $0.563\cdot 10^{-4}$ & $130.00$ \\ \hline
  6400 & $0.221\cdot 10^{-4}$ & $0.142\cdot 10^{-4}$ & $0.220\cdot 10^{-4}$ & $0.142\cdot 10^{-4}$ & $612.60$ \\ \hline
\end{tabular}}
\end{adjustwidth}
\bigskip\bigskip
\end{table}

\subsubsection{Solitary wave run-up}

As the last test case we consider a practically important problem of a solitary wave run-up on a sloping beach. The main difficulty here consists in the fact that the fluid domain may vary in time and its evolution is not known \emph{a priori}. It has to be determined dynamically during the simulation. The wetting/drying algorithm was described in detail in our previous study \cite{Khakimzyanov2016d}. Here we focus mainly on the mesh motion algorithm behaviour.

Consider a closed 1D channel with the bathymetry prescribed by the following function
\begin{equation*}
  z\ =\ -h(x)\ =\ \left\{
  \begin{array}{cl}
    h_0\ -\ x\tan\theta & 0\ \leq\ x\ \leq\ x_s, \\
    -h_0, & x_s\ \leq\ x\ \leq\ \ell,
\end{array}\right.
\end{equation*}
where $\theta$ is the sloping beach angle. The variable bathymetry region is located leftwards from $x_s = 2\cot\theta \approx 38.16$. At the initial time the shoreline is located at $x_{00} = \cot\theta \approx 19.08$.

The initial condition is a solitary wave prescribed by the following formulas (which come from the analytical solitary wave solution to the \textsc{Serre} equations, see \cite{Dutykh2011a} for example):
\begin{align*}
  \eta(x,0)\ &=\ a\;\sech^{2}\biggl[\,\frac{\sqrt{3 g a}}{2 h_0 \sqrt{g(h_0 + a)}}\;(x\ -\ x_{0,w})\biggr], \\ 
  u(x,0)\ &=\ -\sqrt{g(h_0 + a)}\;\frac{\eta(x,0)}{h_0\ +\ \eta(x,0)},
\end{align*}
Variables $a$ and $x_{0,w} = x_s + 30 \approx 68.16$ are the initial solitary wave amplitude and position correspondingly. The right boundary is located at $\ell = x_{0,w} + 30 \approx 98.16$. All numerical parameters used in this experiment are provided in Table~\ref{tab:params5}.

\begin{figure}
  \centering
  \subfigure[]{\includegraphics[width=0.47\textwidth]{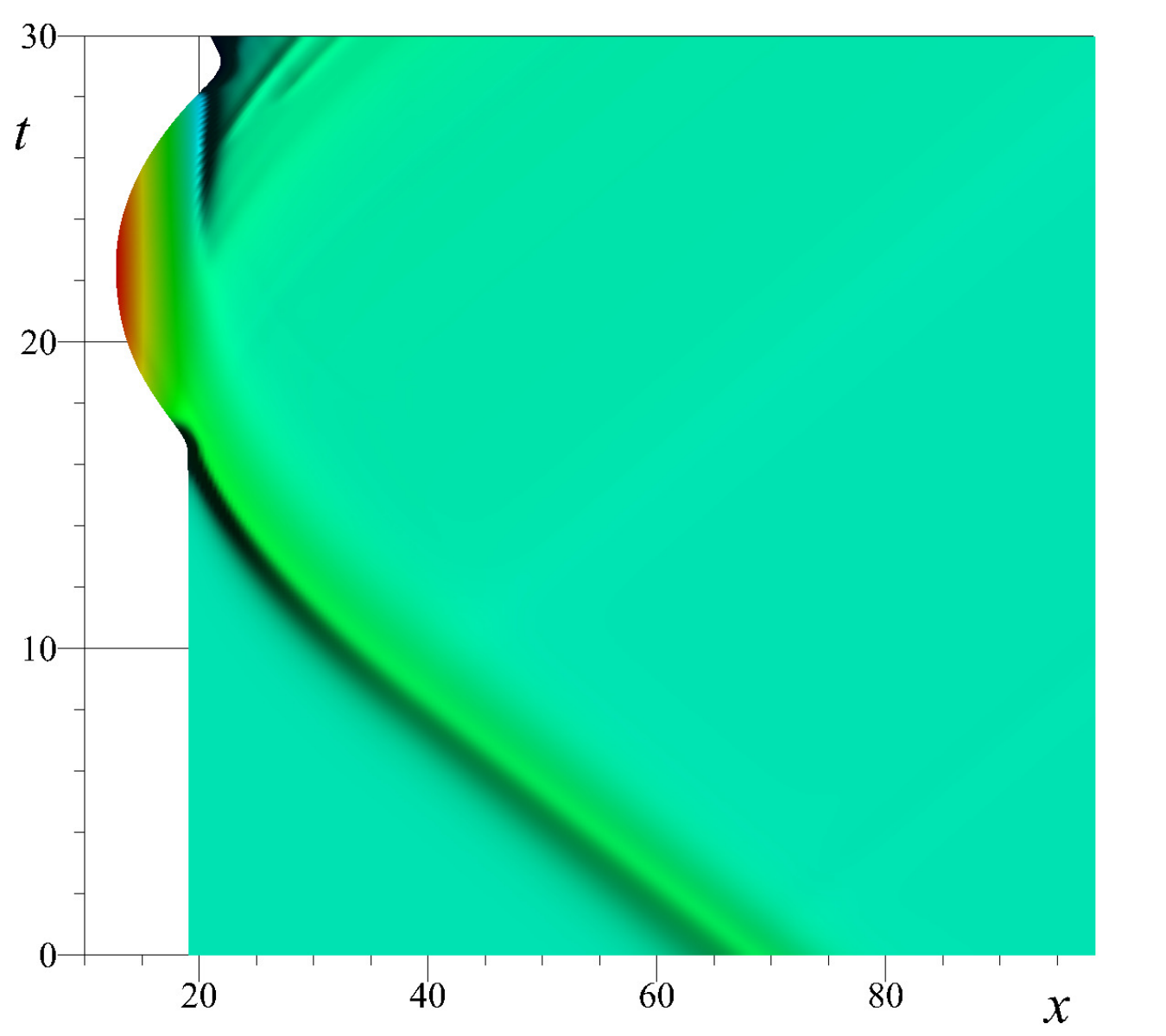}}
  \subfigure[]{\includegraphics[width=0.52\textwidth]{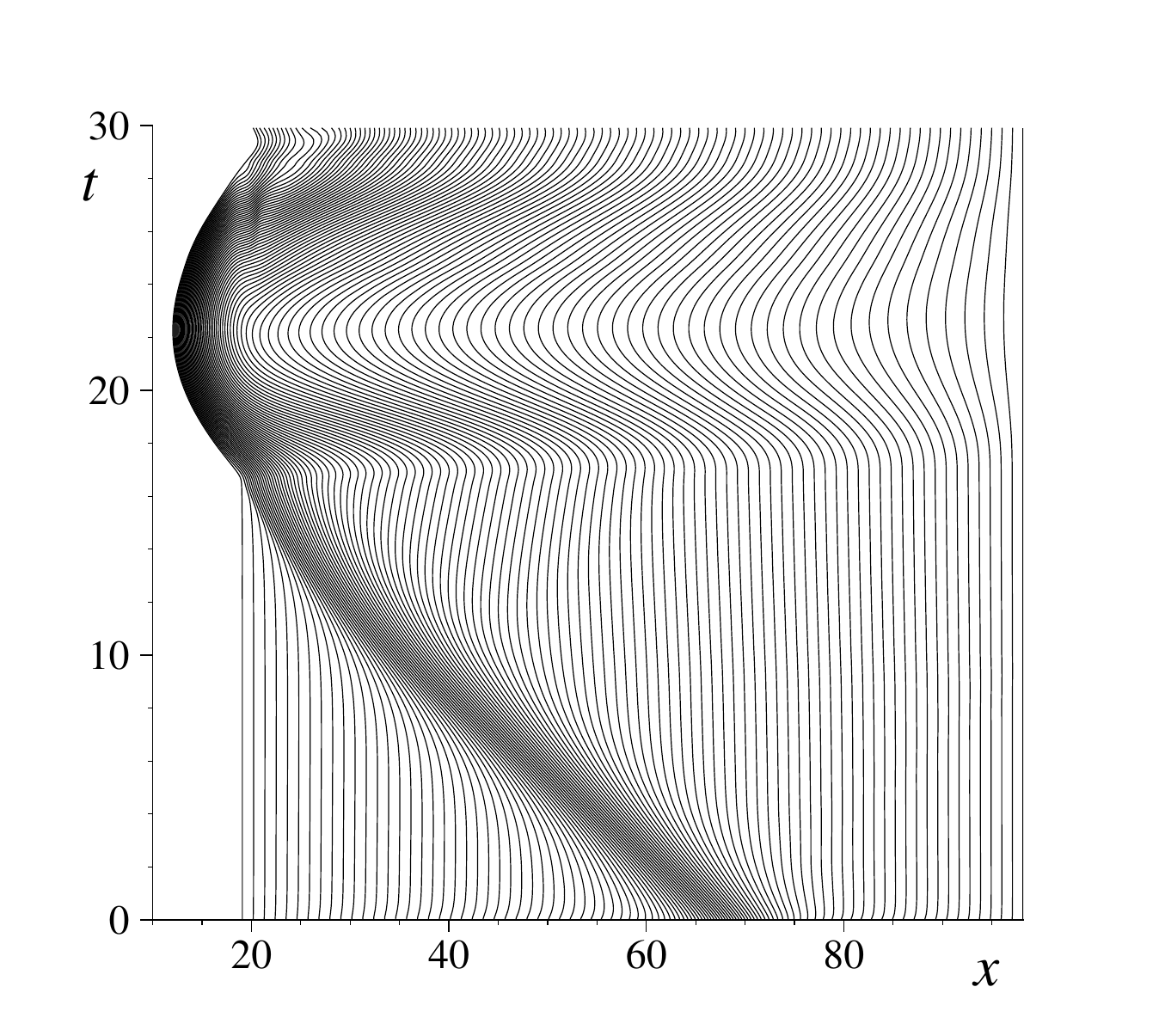}}
  \caption{\small\em Solitary wave run-up simulation on a sloping beach with nonlinear shallow water equations: (a) free surface elevation evolution; (b): trajectories of grid nodes in space-time.}
  \label{fig:6}
\end{figure}

The considered situation is complex enough so that the analytical solution is not available. However, we can simulate it numerically with the predictor--corrector scheme on moving grids presented above. The wave run-up process simulated \emph{in silico} is depicted in Figure~\ref{fig:6}(\textit{a}), while the nodes motion is illustrated in the right panel (\textit{b}). Please, note that only every fifth node is depicted in Figure~\ref{fig:6}(\textit{b}) in order to improve the visibility of the image. The rightmost node is fixed (\ie its trajectory is a straight line), while th leftmost node's trajectory coincides with the shoreline motion. One can see that even on domains with variable spatial extent, the mesh motion algorithm is robust enough and it puts the resolution where it is needed.

The monitor function used in this computation is
\begin{equation}\label{eq:def2}
  w(x,t)\ =\ 1\ +\ \alpha\,\abs{\eta(x,t)}\,.
\end{equation}
It is important to use the wave elevation $\eta(x,t)$ instead of the total water depth $H(x,t)$ because of the solutions such as the lake at rest state. In these situations the total depth $H$ varies in space, while the free surface elevation $\eta(x,t) \equiv 0$. Thus, the usage of the monitor function \eqref{eq:def2} leads to the uniform grid (since $w(x,t)\ \equiv\ 1$) and, trivially, to the well-balanced property shown in Theorem~\ref{thm:well}.

\begin{table}
  \centering
  \begin{tabular}{lc}
  \hline\hline
  \textit{Parameter} & \textit{Value} \\
  \hline
  Computational domain length, $\ell$ & $98.16$ \\
  Number of nodes, $N$ & $1000$ \\
  Final simulation time, $T_f$ & $30.0$ \\
  CFL number, $C$ & $0.95$ \\
  Gravity acceleration, $g$ & $9.81$ \\
  Initial wave amplitude, $a/h_0$ & $0.05$ \\
  Wave crest initial position, $x_{0,w}$ & $68.16$ \\
  Beach slope position, $x_s$ & $38.16$ \\
  Beach slope, $\theta$ & $3^\circ$ \\
  Undisturbed water depth, $h_0$ & $1.0$ \\
  Monitor function parameter, $\alpha$ & $60.0$ \\
  Grid diffusion parameter, $\beta$ & $20.0$ \\
  Smoothing parameter, $\sigma$ & $10.0$ \\
  \hline\hline
  \end{tabular}
  \bigskip
  \caption{\small\em Numerical parameters used in a solitary wave run-up simulation in nonlinear shallow water equations.}
  \label{tab:params5}
\end{table}

Some further validation tests and numerical experiments with this moving grid technique were presented in our preceding work \cite{Khakimzyanov2016d} (without giving many details about the algorithm in use).


\section{Some open problems}
\label{sec:open}

Only in very simple problems (linear and steady) one can study theoretically the choice of the optimal monitor function. Recently we performed such an analysis for a singularly perturbed \textsc{Sturm}--\textsc{Liouville} problem \cite{Khakimzyanov2015b}. In particular, we showed that a judicious choice of the monitor function allows to increase the order of convergence from the 2\up{nd} to the 4\up{th} order (the so-called \emph{supraconvergence} phenomenon), even if the scheme is just 2\up{nd} order accurate on uniform meshes. So, the first open problem consists in providing a similar analysis for, at least, 1D unsteady or nonlinear problem. The method employed in \cite{Khakimzyanov2015b} is suitable for smooth solutions. For discontinuous ones a different technique has to be proposed. Even a linear transport equation is not completely understood from the perspective of constructing a suitable moving grid. It is even possible that for every discretization scheme there is its own ``optimal'' monitor function. So, one choice clearly does not fit all possible schemes. Moreover, even if we choose the monitor function as in \eqref{eq:def2}, it was shown in \cite{Khakimzyanov2015b} that the error may depend non-monotonically on the coefficient $\alpha$. Thus, even if the monitor function form is fixed, one has to set optimal values of coefficients in this function.


\section{Discussion}
\label{sec:concl}

The main conclusions and perspectives of this study are outlined below.

\subsection{Conclusions}

In the present work our goal was to present the moving grid methodology, which does not involve any interpolation procedure (thus, avoiding the unnecessary smearing of the numerical solution). We showed some reasonable choices of the monitor function, which are currently used in practice, but these choices have to be refined further in the forthcoming works. The current state-of-the-art in this field is described in \cite{Huang2011a}.

In this manuscript we presented a detailed description of a family of predictor-corrector schemes up to the second order in space, which were generalized to the case of non-uniform moving meshes, while having all these properties \emph{simultaneously}:
\begin{itemize}
  \item Conservation in space;
  \item Second order accurate (in smooth regions);
  \item TVD property for scalar problems;
  \item Well-balanced character;
  \item Preservation of the discrete conservation law in cell centers and cell interfaces;
  \item Adaptivity in space and in time\footnote{The time step is chosen the largest possible while meeting the stability requirement.};
  \item No interpolation between the meshes, thus, we do not add any unnecessary numerical diffusion;
  \item Implementation is trivial, despite all good properties listed above.
\end{itemize}
This is the main novelty and contribution of our study to the growing field of adaptive numerical methods. This method was illustrated on a number of relatively serious tests coming from linear, nonlinear and vectorial cases. Every time the gain of using a moving grid (with respect to the fixed one) was clearly visible. Consequently, we can only encourage the community to use these (or similar) adaptive strategies to improve the quality and efficiency of their numerical simulations. In particular, the field of tsunami wave modelling could greatly benefit from the moving grid nodes to improve the accuracy and computational efficiency of operational tsunami wave predictions.

\subsection{Perspectives}

The adaptive moving mesh technology presented in this study relies on the choice of the monitor function along with several free real parameters. Numerical experiments show that the method performance may depend crucially on the choice of these real numbers. Consequently, it would be highly desirable to obtain some theoretical estimations of the optimal choice of these parameters at least in the case of the simplest linear advection equation \eqref{eq:adv}. It constitutes one of the main challenges raised by our study. In the absence of such theoretical indications, one might envisage extensive numerical studies in order to issue some practical recommendations on how to choose reasonable values of parameters. Most probably the case of nonlinear systems of conservation laws will remain essentially empiric in the foreseeable future due to the important theoretical difficulties in its analysis.

From a more practical point of view, the predictor-corrector schemes \emph{on moving grids} presented in this study have to be generalized at least to two-dimensional case as it was done for the \textsc{Godunov} scheme in \cite{Azarenok2003}. We already reported an implementation of moving grid technique in two space dimensions for an elliptic problem in \cite{Khakimzyanov2018}. Moreover, the spatial adaptivity has to be coupled with time-adaptive strategies such as the PI step size control \cite{Hairer2009} in order to achieve the full adaptivity.

Finally, there is another goal behind this study. Ultimately, we would like to generalize moving grid methods to some dispersive wave equations \cite{Dutykh2011a}, which are going to be addressed with the operator splitting approach along the lines sketched in \cite{Khakimzyanov2016}. In this way, the hyperbolic part would be addressed with methods outlined above. The particularity of dispersive wave equations is that they possess localised solitary wave solutions. In other words, all the dynamics is concentrated in small portions of space, which move perpetually. Consequently, the grid redistribution techniques seem to be the natural choice to simulate complex solitonic gas dynamics \cite{Dutykh2014d}.


\subsection*{Acknowledgments}
\addcontentsline{toc}{subsection}{Acknowledgments}

This research was supported by RSCF project No 14-17-00219. D.~Dutykh acknowledges the support of the CNRS under the PEPS InPhyNiTi project FARA and project \No EDC26179 --- ``Wave interaction with an obstacle'' as well as the hospitality of the Institute of Computational Technologies of SB RAS during his visit in October 2015. D.~\textsc{Mitsotakis} was supported by the Marsden Fund administered by the Royal Society of \textsc{New Zealand}. We would like to thank Prof. Laurent \textsc{Gosse} (CNR, Italy) for insightful discussions on the topics covered in our study. Moreover, we would like to thank the anonymous Referee on our previous manuscript \cite{Khakimzyanov2016d} who incited us by his comments to produce a study which explains in more details our grid motion technique.


\bigskip\bigskip
\invisiblesection{References}
\bibliographystyle{acm}
\bibliography{biblio}
\bigskip\bigskip


\appendix
\section{Illustrative examples}
\label{app:ex}

Let us consider the `unsteady' version of the equidistribution principle expressed as the following nonlinear elliptic Boundary Value Problem (BVP):
\begin{equation}\label{eq:time}
  \bigl(\,w\,(x,\,t)\;x_{\,q}\,\bigr)_{\,q}\ =\ 0\,, \qquad q\ \in\ ]\,0,\,1\,[\,, \qquad 0\ \leq\ t\ \leq\ T_{\,f}\,,
\end{equation}
with the following boundary conditions:
\begin{equation}\label{eq:bvp}
  x\,(0,\,t)\ =\ 0\,, \qquad x\,(1,\,t)\ =\ \ell\, \qquad 0\ \leq\ t\ \leq\ T_{\,f}\,.
\end{equation}
The first inherent problem with Equation~\eqref{eq:time} is that small perturbations of the monitor function $w\,(x,\,t)$ may lead to finite changes in the solution $x\,(q,\,t)$ as it is illustrated in the following example.


\subsection{Example 1}

Let for $t\ =\ t_{\,n}$ the monitor function $w$ be constant, \ie $w\,(x,\,t_{\,n})\ \equiv\ 1\,$. Then, solution to BVP \eqref{eq:time}, \eqref{eq:bvp} takes the following simple form:
\begin{equation}\label{eq:uni}
  x\,(q,\,t_{\,n})\ =\ q\,\ell\,, \qquad q\ \in\ [\,0,\,1\,]\,.
\end{equation}
This solution prescribes the uniform distribution of nodes on $[\,0,\,\ell\,]\,$. Let us also assume that at the next time step $t\ =\ t_{\,n+1}$ the monitor function $w$ was perturbed with quantity $\eps$ on a portion of the segment $[\,0,\,\ell\,]\,$:
\begin{equation*}
  w\,(x,\,t_{\,n+1})\ =\ \begin{dcases}
    \ 1\ +\ \eps\,, \quad & x\ \leq\ x_{\,\ast}\,, \\
    \ 1\,, \quad & x\ >\ x_{\,\ast}\,,
  \end{dcases}
  \qquad x_{\,\ast}\ \in\ ]\,0,\,\ell\,[\,.
\end{equation*}
Then, solution to BVP \eqref{eq:time}, \eqref{eq:bvp} will be given by the following piece-wise linear function:
\begin{equation}\label{eq:sol}
  x\,(q,\,t_{\,n+1})\ =\ \begin{dcases}
    \ q\;\frac{x_{\,\ast}}{q_{\,\ast}}\,, \quad & q\ \leq\ q_{\,\ast}\,, \\
    \ x_{\,\ast}\ +\ (q\ -\ q_{\,\ast})\;\frac{\ell\ -\ x_{\,\ast}}{1\ -\ q_{\,\ast}}\,, \quad & q\ >\ q_{\,\ast}\,,
  \end{dcases}
\end{equation}
where $q_{\,\ast}\ =\ \dfrac{(1\ +\ \eps)\,x_{\,\ast}}{\ell\ +\ \eps\,x_{\,\ast}}\,$. From here, the following estimation follows:
\begin{equation*}
  \norm{x(\cdot,\,t_{\,n+1})\ -\ x(\cdot,\,t_{\,n})}_{\,C\,(Q)}\ =\ \abs{x(q_{\,\ast},\,t_{\,n+1})\ -\ x(q_{\,\ast},\,t_{\,n})}\ =\ \frac{\eps\,(\ell\ -\ x_{\,\ast})\,x_{\,\ast}}{\ell\ +\ \eps\,x_{\,\ast}}\,.
\end{equation*}
It is not difficult to check that this norm of the difference takes the maximal value for $x_{\,\ast}^{\,\max}\ =\ \dfrac{\ell\,(\sqrt{1\ +\ \eps}\ -\ 1)}{\eps}\,$. For this given $x_{\,\ast}^{\,\max}$ and for small values of $\eps$ we have:
\begin{equation*}
  \norm{x(\cdot,\,t_{\,n+1})\ -\ x(\cdot,\,t_{\,n})}_{\,C\,(Q)}\ =\ \frac{\ell\,(2\ +\ \eps\ -\ 2\,\sqrt{1\ +\ \eps})}{\eps}\ =\ \frac{\eps\,\ell}{4}\ +\ \O\,(\eps^{\,2})\,.
\end{equation*}
Hence, by choosing the domain length $\ell$ sufficiently large, we can obtain a finite\footnote{To obtain truly finite perturbations in $x\,(q,\,t)\,$, it is sufficient to take $\ell\ =\ \O\,(\eps^{\,-1})\,$.} perturbation in the solution to BVP \eqref{eq:time}, \eqref{eq:bvp}.

In numerical practice, it means that the grid nodes in the vicinity of the point $x\ =\ x_{\,\ast}^{\,\max}$ will undergo a shift proportional to $\propto\ \eps\,\ell$ when we advance the solution from the time layer $t\ =\ t_{\,n}$ to $t\ =\ t_{\,n+1}\,$. If, for some reason, on the time layer $t\ =\ t_{\,n+1}$ the monitor function $w$ becomes constant\footnote{In this case, the monitor function may assume two possible constant values $w\,(x,\,t)\ \equiv\ 1$ or $w\,(x,\,t)\ \equiv\ 1\ +\ \eps\,$.} again, the grid nodes will experience an opposite move towards a uniform grid \eqref{eq:uni} again. Such oscillations have been observed numerically and the present example provides a theoretical explanation of this phenomenon.


\subsection{Example 2}

In order to complete the description of possible shortcomings of the equidistribution principle for the mesh motion in evolutionary PDEs, in this Section we consider a simple unsteady problem where the scenario described in Example~1 is actually realized.

Consider the linear advection equation:
\begin{equation*}
  u_{\,t}\ +\ a\,u_{\,x}\ =\ f\,(x,\,t)\,, \qquad a\ =\ \const\ >\ 0\,,
\end{equation*}
where the source term is chosen to be
\begin{equation*}
  f\,(x,\,t)\ =\ \begin{dcases}
  \ \eps\,(x_{\,\ast}\ -\ x\ -\ a\,t)\,, \quad & x\ \leq\ x_{\,\ast}\,, \\
  \ 0\, \quad & x\ >\ x_{\,\ast}\,.
  \end{dcases}
\end{equation*}
In this Section $\eps$ is a small parameter and $x_{\,\ast}\ \in\ ]\,0,\,\ell\,[\,$. An exact solution to the homogeneous Initial-Boundary Value Problem (IBVP) in domain $[\,0,\,\ell\,]\,\setminus\,\{\,x_{\,\ast}\,\}$ is given by the following formula:
\begin{equation*}
  u\,(x,\,t)\ =\ \begin{dcases}
  \ \eps\,t\,(x_{\,\ast}\ -\ x)\,, \quad & x\ <\ x_{\,\ast}\,, \\
  \ 0\,, \quad & x\ >\ x_{\,\ast}\,.
  \end{dcases}
\end{equation*}
Let us choose the monitor function as
\begin{equation*}
  w\,(x,\,t)\ =\ 1\ +\ \alpha\,\abs{u_{\,x}\,(x,\,t)}
\end{equation*}
and take the parameter $\alpha\ \equiv\ 1$ for simplicity. Then, it is not difficult to check that
\begin{equation*}
  w\,(x,\,t)\ =\ \begin{dcases}
  \ 1\ +\ \eps\,t\,, \quad & x\ <\ x_{\,\ast}\,, \\
  \ 1\,, \quad & x\ >\ x_{\,\ast}\,.
  \end{dcases}
\end{equation*}
The non-uniform grid given by the equidistribution principle is obtained as a solution to BVP \eqref{eq:time}, \eqref{eq:bvp}. In this example, it takes precisely the form \eqref{eq:sol}, but the point $q_{\,\ast}$ depends on time:
\begin{equation*}
  q_{\,\ast}\ =\ \frac{(1\ +\ \eps\,t)\,x_{\,\ast}}{\ell\ +\ \eps\,t\,x_{\,\ast}}\,.
\end{equation*}
It is not difficult to show that the maximal speed of grid nodes $-\,x_{\,t}$ is achieved at $t\ =\ 0$ and $x_{\,\ast}\ =\ \frac{\ell}{2}\,$. Moreover, $\max\{-\,x_{\,t}\}\ =\ \frac{\eps\,\ell}{4}\,$. The grid at $t\ =\ 0$ is uniform and at the first time step it will experience a brutal shift to the left. The maximal displacement being equal to $\dfrac{\eps\,\ell}{4}\,$. This shift may provoke, at least, a loss of accuracy or, in the worst case scenario, blow-up of the numerical solution. The latter may routinely happen in some explicit schemes, where the time step does not take into account the grid nodes displacements.


Henceforth, we have just demonstrated analytically that the grid constructed with equidistribution principle~\eqref{eq:time} turns out to be unstable when the equidistribution principle is applied in unsteady problems. Indeed, arbitrary small perturbations in the monitor function or in the governing equation may lead to finite perturbations in the grid. Thus, we highly recommend to use the unsteady parabolic Equation~\eqref{eq:parab} in order to move grid nodes smoothly.


\section{Abbreviations}

\begin{description}
  \item[1D] One-dimensional
  \item[2D] Two-dimensional
  \item[BVP] Boundary Value Problem
  \item[CFL] \textsc{Courant}--\textsc{Friedrichs}--\textsc{Lewy}
  \item[DART] Deep-ocean Assessment and Reporting Tsunamis
  \item[FEM] Finite Elements Method
  \item[FNWD] Fully Nonlinear Weakly Dispersive
  \item[IBVP] Initial-Boundary Value Problem
  \item[MOST] Method Of Splitting Tsunami
  \item[NSWE] Nonlinear Shallow Water Equations
  \item[PDE] Partial Differential Equation
  \item[TVD] Total Variation Diminishing
  \item[WNWD] Weakly Nonlinear Weakly Dispersive
\end{description}

\end{document}